%% file: Improved_Catoni_CS.tex
\documentclass[reqno,final,11pt]{alea-amsart}
\usepackage{natbib}  
\usepackage{fancyhdr} 
\usepackage{color} 
\usepackage{hyperref} 
\usepackage{graphicx} 
\usepackage{booktabs}
\usepackage{multirow}
\usepackage{tipa}
\usepackage{mathrsfs}  
\usepackage{amsfonts}
\usepackage{amssymb}
\usepackage{amsmath}
\usepackage{extarrows}
\usepackage{amsthm}
\usepackage{etoolbox}
\numberwithin{equation}{section}
\usepackage{amsbsy}
\usepackage{subcaption}
\usepackage{graphics}
\usepackage{placeins}  
\usepackage{setspace}
\usepackage{natbib} 
\usepackage{enumitem} 
\usepackage{float}
\usepackage{caption}
\usepackage{geometry}
\usepackage[utf8]{inputenc}
\usepackage[T1]{fontenc}
\usepackage{verbatim}

\usepackage{tikz}
\usetikzlibrary{calc,shapes,backgrounds,arrows}
\usepackage{amssymb}

\definecolor{aleacolor}{rgb}{0.16,0.59,0.78}

\hypersetup{
breaklinks,
colorlinks=true,
linkcolor=aleacolor,
urlcolor=aleacolor,
citecolor=aleacolor}

\renewcommand{\cite}{\citet}

\theoremstyle{plain}
\newtheorem{theorem}{Theorem}[section]                                          
                          
\newtheorem{lemma}[theorem]{Lemma}
\newtheorem{corollary}[theorem]{Corollary}

\theoremstyle{definition}
\newtheorem{definition}[theorem]{Definition}
\theoremstyle{remark}
\newtheorem{remark}[theorem]{Remark}
\newtheorem{example}[theorem]{Example}

\renewcommand{\qedsymbol}{\ensuremath{\Box}}
\AtBeginEnvironment{definition}{\pushQED{\qed}}
\AtEndEnvironment{definition}{\popQED}
\AtBeginEnvironment{remark}{\pushQED{\qed}}
\AtEndEnvironment{remark}{\popQED}
\AtBeginEnvironment{example}{\pushQED{\qed}}
\AtEndEnvironment{example}{\popQED}

\makeatletter
\let\catoni@original@equation@qed\equation@qed
\def\equation@qed{%
  \iftagsleft@
    \catoni@original@equation@qed
  \else
    \gdef\alt@tag{%
      \global\let\alt@tag\@empty
      \hbox{\tagform@\theequation\quad\qedsymbol}%
      \setbox\z@
    }%
  \fi
}
\makeatother

\makeatletter \@addtoreset{equation}{section} \makeatother
\renewcommand\theequation{\thesection.\arabic{equation}}
\numberwithin{figure}{section}
\numberwithin{table}{section}

\renewcommand{\ejournal}{ALEA, Lat. Am. J. Probab. Math. Stat.}
\makeatletter
\renewcommand{\@seteheader}{%
  \ \vskip -1cm
  \noindent\ejournal\hfill\raisebox{0pt}[39.05048pt][33.5755pt]{}\\
  \vskip -1.2cm
  \noindent\raisebox{0pt}[7.54317pt][2.12914pt]{}%
  \vskip 1cm
}
\makeatother

\usepackage{xcolor}
\DeclareRobustCommand{\revision}[1]{\ifvmode\leavevmode\fi#1}

\makeatletter
\expandafter\def\csname revisionbib@FangCS\endcsname{1}
\expandafter\def\csname revisionbib@WangSGD2021\endcsname{1}
\AtBeginEnvironment{thebibliography}{%
  \let\revision@original@bibitem\bibitem
  \renewcommand{\bibitem}[2][]{%
    \revision@original@bibitem[#1]{#2}%
    \leavevmode
    \ifcsname revisionbib@#2\endcsname\color{black}\else\color{black}\fi\ignorespaces}}
\makeatother
\makeatletter
\expandafter\def\csname revisioneq@eq:stitched-admissible-times\endcsname{1}
\expandafter\def\csname revisioneq@eq:stitched-error-budget\endcsname{1}
\expandafter\def\csname revisioneq@eq:stitched-failure-events\endcsname{1}
\expandafter\def\csname revisioneq@eq:stitched-threshold\endcsname{1}
\expandafter\def\csname revisioneq@eq:width-envelopes\endcsname{1}
\expandafter\def\csname revisioneq@eq:width-root-coefficients\endcsname{1}
\expandafter\def\csname revisioneq@eq:width-score\endcsname{1}
\let\revision@original@tagform\tagform@
\DeclareRobustCommand{\eqref}[1]{\textup{\revision@original@tagform{\ref{#1}}}}
\def\tagform@#1{%
  \ifcsname revisioneq@\df@label\endcsname
    \maketag@@@{\revision{(\ignorespaces#1\unskip\@@italiccorr)}}%
  \else
    \revision@original@tagform{#1}%
  \fi}
\makeatother
\begin{document}

\title[Improved Catoni-type confidence sequences with infinite variance]{Improved Catoni-type confidence sequences for \\ estimating the mean when the variance is infinite}

\author{Chengfu Wei}
\address{Institute for Financial Studies, Shandong University, Jinan 250100, China}
\email{chengfuwei1024@gmail.com}

\author{Jordan Stoyanov}
\address{Institute of Mathematics $\&$ Informatics, Bulgarian Academy of Sciences, 1113 Sofia, Bulgaria, \newline and~Taishan College, Shandong University, Jinan 250100, China}
\email{stoyanovj@gmail.com}

\author{Yiming Chen}
\address{School of Mathematical Sciences, Peking University, Beijing 100871, China}
\email{ymchenmath@math.pku.edu.cn}

\author{Zijun Chen}
\address{Institute for Financial Studies, Shandong University, Jinan 250100, China}
\email{czj4096@gmail.com}

\thanks{Chen Yiming was  supported by  the National Natural Science Foundation of China (Grant no. 12071257).}
\subjclass[2010]{60G42, 62F25, 62L12} 
\keywords{Heavy-tailed distribution,  Estimating the mean,  Supermartingales, Confidence sequence, Improved Catoni-type estimator, Stitching method}

\begin{abstract}
  We construct an improved Catoni-type confidence sequence (CS) for estimating the unknown constant conditional mean of an adapted discrete-time stochastic process whose variance may be infinite but whose conditional absolute central moment of order $(1+\alpha)$, $\alpha\in(0,1]$, has a known uniform bound. Our model is beyond the traditional i.i.d. case. The construction is developed within the Catoni-type supermartingale framework and incorporates a sharper calibration of the influence function.

  Under stated sufficient conditions, we obtain improved explicit width bounds. The standard stitching method yields a CS with coverage from the first observation and a high-probability time-uniform width bound beyond a computable deterministic threshold. For fixed error level $\delta$, the bound has a shrinkage rate of $\mathcal{O}\big(((\log\log t)/t)^{\alpha/(1+\alpha)}\big)$ as the sample size $t\to\infty$, and the width of the same CS converges to zero almost surely at this rate. Within the stated admissible range, the bound has a growth rate of $\mathcal{O}\big((\log(1/\delta))^{\alpha/(1+\alpha)}\big)$ as the error level $\delta\to0$. At the benchmark calibration and for $\alpha\in(0,1)$, the proposed bound has a smaller width constant than the comparison bound without increasing the sufficient sample-size threshold. Simulations under two heavy-tailed models with infinite variance show smaller empirical mean interval widths for the proposed CSs.
\end{abstract}

\maketitle

\input{section/introduction}

\input{section/Problemset-up}

\input{section/Mainresults}

\input{section/Simulationresults}

\appendix
\input{section/Appendix}

\section*{Acknowledgements}
Three of the authors, CW, YC and ZC, would like to express their sincere gratitude to their supervisor, Professor Hanchao Wang, for his invaluable guidance and continuing encouragement.

One of the authors, JS, acknowledges the warm hospitality and support provided by Taishan College, Shandong University, during his visits in 2024 and 2025.

We are also indebted to Professor Luc Devroye for his insightful comments on the first version of the manuscript after it was posted online; see \href{https://arxiv.org/abs/2409.04198}{arXiv:2409.04198}.

Finally, we are grateful to the anonymous reviewers and the editors of \textit{ALEA} for their detailed constructive comments and suggestions, which helped us to improve the clarity, enhance the content and  make the whole text read more smoothly. In particular,   the reviewer of the revised manuscript  drew  our attention to the forthcoming paper \cite{FangCS} and advised 
 us to clarify the similarities and differences between our findings and the results in that paper. These points have been addressed in the revised manuscript. 

The authors take full responsibility for the content of this article, including all mathematical statements, proofs, numerical results, and figures.

\bibliographystyle{alea3}
\bibliography{reference}

\end{document}

%% file: section/introduction.tex
\section{Introduction}\label{introduction}
Inference on the population mean from finite samples is a fundamental problem in applied probability, statistics and machine learning.
The present paper studies confidence sequences (CSs) for heavy-tailed mean estimation, with particular emphasis on settings in which the variance may be infinite but the central moment of order $(1+\alpha)$ admits a known upper bound.
We first briefly review some developments of CSs for mean estimation and their extension to heavy-tailed models. Then we discuss how our work relates to the broader literature on robust mean estimation, and finally summarize our main contribution.

For convenience, we recall some common notations that will be used throughout the paper. Let $\mathbb{N}=\{1,2,\ldots\}$ and $\mathbb{N}_0=\{0,1,2,\ldots\}$ denote the sets of positive and nonnegative integers, respectively. For a nonempty interval $I\subseteq\mathbb{R}$, its lower and upper endpoints are $\inf I$ and $\sup I$, respectively, allowing extended real values. We use the conventions $\inf\varnothing=+\infty$ and $\sup\varnothing=-\infty$. The width $|I|$ of the interval $I$ is the difference $\sup I-\inf I$ for a nonempty bounded interval; it is $0$ for the empty set, and $+\infty$ for an unbounded interval. The notation `$\log$' denotes the natural logarithm. We also use the standard asymptotic conventions: for two positive real number sequences $\{x_t\}_{t=1}^{\infty}$ and $\{y_t\}_{t=1}^{\infty}$, we write:

$x_t=\mathcal{O}(y_t)$ \ if \ $\limsup_{t\to\infty} x_t/y_t<\infty$;

$x_t=\widetilde{\mathcal{O}}(y_t)$ \ if \ both $x_t=\mathcal{O}(y_t)$ and $y_t=\mathcal{O}(x_t)$ hold;

$x_t=o(y_t)$ \ if \ $\lim_{t\to\infty}x_t/y_t=0$;

$x_t \asymp y_t$ \ if \ $\lim_{t\to\infty} x_t/y_t$ exists and is positive and finite.

\smallskip
Theorems, Lemmas, Corollaries, Definitions, Remarks, and Examples share a counter and are numbered consecutively within each section, with the section number preceding the item number.

The symbol \ $\Box$ \  marks the end of each Definition, Remark, Example, and Proof.

\subsection{Confidence sequences and heavy-tailed mean estimation}\label{CS}~\\

We begin with  the classical fixed-time problem of mean estimation. Let $X$ be a real-valued random variable defined on an underlying probability space $(\Omega,\mathcal{F},\mathsf{P})$, and write $X\sim F$, where $F$ denotes the distribution function of $X$. For any $t\in\mathbb{N}$, consider i.i.d.\ observations $X_1,X_2,\ldots,X_t$ drawn from $F$. The goal is to estimate the unknown mean $\mu$ of $X$:
\begin{equation*}
    \mu = \int_{\Omega} X(\omega)\mathrm{~d}\mathsf{P}(\omega)
    = \int_{\mathbb{R}} x \, \mathrm{d} F(x)
    = \mathsf{E}[X]
    \quad (= \mathsf{E}[X_t] \ \mbox{for any } t).
\end{equation*}

A classical approach to this problem is to construct a confidence interval (CI) based on the empirical mean computed from the available observations.

In the sequel, we use the notation $\mathcal{F}_t := \sigma(X_1, \ldots, X_t)$ for the sigma-algebra generated by the first $t$ observations, and say that $\{\mathcal{F}_t\}_{t=1}^{\infty}$ is the natural filtration associated with the observation sequence $X_1, X_2, \ldots$.

For a given time $t$, this procedure requires two steps:
\begin{enumerate}[label=(\roman*), nosep, leftmargin=1.5em, labelsep=0.5em]
    \item Define an estimator $\widehat{\mu}_t$ of $\mu$, where $\widehat{\mu}_t = \widehat{\mu}_t(X_1,\ldots, X_t)$ is $\mathcal{F}_t$-measurable;
    \item Construct a confidence interval $\mathrm{CI}_t$ containing $\mu$ with confidence level $1-\delta$, where $\delta \in (0,1)$ is the error level.
\end{enumerate}

\begin{definition}[Confidence interval, CI] \label{CI definition}
    For any $t \in \mathbb{N}$ and fixed $\delta \in (0,1)$, the $\mathcal{F}_t$-adapted random set $\mathrm{CI}_t$ on $\mathbb{R}$ is called a $(1-\delta)$-confidence interval for the unknown mean $\mu$ if it satisfies
    \begin{equation}\label{CI}
        \mathsf{P}[\,\omega: \mu \in \mathrm{CI}_{t}(\omega)] \ \geqslant \ 1-\delta, \quad t \in \mathbb{N}.\qedhere
    \end{equation}
\end{definition}

In general, $\mathrm{CI}_t$ is a random interval depending on the sample $X_1,\ldots, X_t$. Of principal importance is its {\bf \emph{width}}, denoted by $|\mathrm{CI}_t|$. The width is nonnegative and depends on both the sample size $t$ and the error level $\delta$.

In both theory and applications there are variations of Definition \ref{CI definition}. For instance, in sequential testing problems, instead of \eqref{CI}, one considers relations of the form
\[
    \mathsf{P}[\,\omega: \mu \in \mathrm{CI}_{\tau(\omega)}(\omega)] \  \geqslant \ 1-\delta,
\]
where $\tau = \tau (\omega)$ is a stopping time with respect to the filtration $\mathcal{F}_t$. Some delicate analytic questions arise. These can be addressed by adopting a stronger definition than Definition \ref{CI definition}:

\begin{definition}[Confidence sequence, CS] \label{CS definition}
    For fixed $\delta \in (0,1)$, the sequence of random sets $\{\mathrm{CI}_t\}_{t=1}^{\infty}$ on $\mathbb{R}$, where each $\mathrm{CI}_t$ is $\mathcal{F}_t$-measurable, is called a $(1 - \delta)$-confidence sequence for the unknown mean $\mu$ if
    \begin{equation*}
        \mathsf{P}[\,\omega: \mu \in \mathrm{CI}_{t}(\omega) \ \mbox{for all } \ t \in \mathbb{N}\,] \ \geqslant \ 1 - \delta.\qedhere
    \end{equation*}
\end{definition}

\begin{remark}[Sequential coverage and interval widths]
    Definition \ref{CS definition} permits the confidence bounds to be recalculated and updated repeatedly as new observations become available, while maintaining valid statistical coverage throughout time. For more detailed discussions of sequential testing and CSs, see \cite{HR2021} and \cite{JK2017}. Throughout the paper, when  referring to the width of a CS, we mean the widths $|\mathrm{CI}_t|, t \in \mathbb{N}$ of the random intervals constituting the whole sequence.
\end{remark}

The primary focus of our study is on the behavior of the width $|\mathrm{CI}_t|$ as a function of the sample size $t$ and the error level $\delta$. Specifically, we examine the following two aspects:
\begin{enumerate}[label=(\alph*), nosep, leftmargin=1.5em, labelsep=0.5em]
    \item The rate of shrinkage describing how quickly $|\mathrm{CI}_t|$ decreases as $t \to \infty$;
    \item The rate of growth describing how quickly $|\mathrm{CI}_t|$ increases as $\delta \to 0$.
\end{enumerate}

Previous studies on CSs for the unknown mean $\mu$ typically relied on rather restrictive assumptions on the underlying distribution $F$. Early works such as \cite{DR1967} and \cite{JT1989} focused on Gaussian settings, while \cite{K1978} extended this line to the multivariate normal case. A notable step forward was made by \cite{L1976a}, who constructed CSs for several common parametric families, including Bernoulli and exponential distributions, and later generalized the approach to exponential families of distributions. See also \cite{CA1979} and \cite{J2015} for subsequent extensions.

In parallel, \cite{L1976b} introduced the conceptually important ``power-one'' approach to sequential testing, which led to a foundational construction of CSs for light-tailed distributions. More specifically, that work assumed the existence of the moment generating function (m.g.f.), defined formally as $M(t) := \mathsf{E}\,[\mathrm{e}^{tX}] < \infty$ for all $t \in (-t_0, t_0)$ with some $t_0>0$. This assumption, known as Cram\'er's condition, characterizes light-tailed distributions and ensures that all moments exist and uniquely determine the distribution; see \cite{JS2012} for further discussion.

More recently, \cite{HR2020,HR2021} substantially broadened the scope of CSs through a systematic nonparametric development. In place of earlier moment-generating-function-based assumptions, they introduced a general sub-$\psi$ framework, which covers many common light-tailed classes of distributions, including the sub-Bernoulli, sub-Gaussian, sub-Poisson, sub-Gamma, and sub-exponential families.

In contrast to the above literature, which largely relies on light-tail assumptions such as sub-Gaussianity, \cite[Theorem 9]{CS} considerably weakened the distributional requirements. They constructed a Catoni-type CS based on the Catoni {\rm M}-estimator of \cite{Catoni}, assuming only a known upper bound on the variance $\sigma^2=\mathsf{Var}[X]$. Even under this minimal finite-variance assumption, their proposed Catoni-type CSs achieve significant width control. More recently, \cite{Huber} further developed robust Catoni-style CSs under Huber's contamination model, in which it is assumed that the observations may contain an $\varepsilon$-fraction of arbitrary corruptions and that only a known upper bound on the $(1+\alpha)$-th central moment of the target distribution is available.

The importance of the Catoni-type CSs of \cite{CS} in the finite-variance setting is that replacing the sub-Gaussian assumption by the weaker assumption of merely a finite variance does not change the order of the  width. Indeed, these authors achieved CSs of the same growth rate of order $\mathcal{O}(\sqrt{\log(1/\delta)})$, and the same shrinkage rate of order $\mathcal{O}(\sqrt{(\log\log t)/t})$, as the sub-Gaussian CSs in \cite{HR2021}.  Notice that the latter rate matches the classical law-of-the-iterated-logarithm lower bound. Moreover, the experimental results in \cite{CS} indicate that the empirical performance of sub-Gaussian CSs can be `close' to that of Catoni-type CSs. At the same time, however, sub-Gaussian assumptions are often difficult to verify in practice and clearly fail for several heavy-tailed models.  This observation further motivates the study of robust CSs under weaker moment conditions.

Building on the fixed-time work of \cite{Chen}, \cite[Theorem 12]{CS} further extended the Catoni-type approach to the setting in which the variance may be infinite, under the assumption that  for some $\alpha\in(0,1)$, a known upper bound on the central $(1+\alpha)$-th moment is available. In this way, they constructed an $\alpha$-Catoni-type CS. Such heavy-tailed models are well motivated and arise naturally in a wide range of applications; see \cite{RS2007} for further discussion.

In addition to the sequential literature discussed above, the present work is also closely related to studies on fixed-time Catoni-type mean estimation in the infinite-variance setting, especially to the refinement in \cite{NOC}, which is built on \cite{Chen}. In particular, \cite{Chen} proposed an estimator whose rate of growth is close to that of the Catoni {\rm M}-estimator as $\alpha\to1$. However, these authors assumed a sub-optimal influence function coming from Taylor's expansion, employed $C_r$-inequalities, and imprecisely characterized the roots of certain equations involving powers of order less than $2$ to construct their estimator. These loose intermediate steps led to larger asymptotic constants than those in \cite{Catoni}, thus precluding near-optimality. Through a different approach, \cite{NOC} mitigated these shortcomings and proposed an estimator for $\mu$ when $\alpha<1$ that nearly optimizes the Catoni {\rm M}-estimator. The proposed estimator has the same order and asymptotic constants as the Catoni {\rm M}-estimator in \cite{Catoni}.

\subsection{Relation to robust mean estimation and sequential inference}\label{robust estimation}~\\

Available is a rich literature on robust mean estimation under heavy tails with several classical lines of work including Catoni {\rm M}-estimators, median-of-means (MoM) estimators, trimmed means, and estimators based on truncation or thresholding. In particular, \cite{Catoni} introduced a PAC-Bayesian {\rm M}-estimation approach with sub-Gaussian-type deviation bounds, which is a   finite-variance model. MoM-type procedures were subsequently developed into a major paradigm for heavy-tailed mean estimation; see the survey of \cite{LM2019}. In parallel, \cite{CG2018} proposed a dimension-free thresholding estimator for random vectors, while \cite{LM2021trim} established the optimality of trimmed mean in a robust multivariate setting. More recently, \cite{OOR2025} proved sharp finite-sample properties of the trimmed mean under finite variance, showing a sub-Gaussian-type behavior in the scalar case.

At a generic level, fixed-time deviation bounds can often be converted into time-uniform estimation procedures via general wrappers, typically at the cost of an additional iterated-logarithm factor. One example is the doubling-based conversion in \cite{DH2024}, which is used as a comparison device in \cite{CWR2025}. However, this type of conversion to time-uniform inference is essentially i.i.d.-based: the construction in \cite{DH2024} is formulated for independent and identically distributed observations, and \cite{CWR2025} likewise explicitly acknowledged this limitation when using that method as a benchmark.
Thus, such a generic time-uniform conversion should be distinguished from a direct CS theory for a specific robust estimator family under more general dependence assumptions, such as martingale dependence. With this distinction in mind, it is natural to separate the existing literature into methods for which a direct CS theory is already available, and methods that currently enter the sequential setting mainly through indirect i.i.d.-based conversions of fixed-time results.

Within the first category, two robust routes are especially relevant here. The first is the Catoni-style {\rm M}-estimation route of \cite{CS,Huber}, which constructs CSs for a scalar heavy-tailed mean via robust exponential supermartingales, and is therefore naturally compatible with Ville's inequality and stitching arguments. The second is the closely related but technically distinct Catoni--Giulini truncation route. In this direction, \cite{CWR2025} derived time-uniform confidence sphere sequences for random vectors and, in the heavy-tailed setting, provided sequential extensions of the Catoni--Giulini estimator under finite-variance assumptions. More recently, \cite{WCMR2025} further extended this truncation-based line to $2$-smooth Banach spaces under central or non-central $(1+\alpha)$-th moment assumptions for some $\alpha \in(0,1]$.
This allows one to deal with infinite variance and martingale dependence, and to establish both line-crossing inequalities and non-asymptotic law-of-the-iterated-logarithm bounds. Thus, although the Catoni {\rm M}-estimator line and the Catoni--Giulini truncation line look generally close, they are technically distinct: the former is influence-function-based, whereas the latter is truncation-based.

By contrast, MoM and trimmed mean currently fall mainly into the second category, entering the sequential literature through such indirect conversions rather than through direct CS constructions. As discussed in \cite{CWR2025}, the MoM baselines are made time-uniform through generic wrappers such as naive union bounds or the doubling method, rather than through a direct construction of CSs based on MoM. Their experiments suggest that the stitched Catoni--Giulini-type procedure performs best overall, although the  geometric MoM remains competitive in some regimes. Likewise, for trimmed means, direct non-asymptotic CS constructions appear to be less developed in the current literature. Under Huber's contamination model, \cite{Huber} compared directly constructed robust CSs with sequential baselines induced from the fixed-time trimmed-mean route of \cite{LM2021trim}, and showed that the direct Catoni-style supermartingale construction can yield substantially tighter intervals.

During the revision of this paper, the reviewer drew our attention to the forthcoming journal version of \cite{FangCS}. These authors study Catoni-type CSs under a known bound on the $(1+\alpha)$th absolute central moment, $0<\alpha\leqslant1$. Their main results are formulated and proved for i.i.d. observations, while their Remark~1 states that their results `carry over for stochastic processes with constant conditional expectation using the influence function and the standard supermartingale arguments'. They refine the width analysis using Ville's inequality, establish lower bounds for specified Catoni-type constructions, and obtain time-uniform stitched width bounds of iterated-logarithmic order. Both works use the influence-function calibration of \cite{NOC}.

In the present paper, we study a quite general martingale-type model with constant conditional mean and bounded conditional $(1+\alpha)$th absolute central moment. Under assumptions (A1)--(A2) in Section~\ref{problem set-up}, we perform a detailed analysis of two-sided non-asymptotic CSs, with an emphasis on explicit width constants and computable sufficient thresholds. At the benchmark parameter choices and for $0<\alpha<1$, these explicit width upper bounds are sharper than the corresponding comparison bounds for the constructions of \cite{CS}, whenever the respective sufficient conditions hold. The correspondence with the supermartingale constructions of \cite{FangCS} and the specialization of their fixed-time width analysis are discussed before Lemma~\ref{catoni supermartingales} and after Theorem~\ref{width of CS}, respectively.

\subsection{Our contributions}\label{contributions}~\\

Building on the Catoni-type sequential framework of \cite{CS} and the standard stitching construction of \cite{HR2021}, the present paper makes several contributions in the heavy-tailed setting with possibly infinite variance.

$\Diamond$ \ The model we deal with is quite general and beyond the traditional i.i.d. case. Following the conditional-mean framework of \cite{CS}, we consider  adapted discrete-time stochastic processes of martingale type, whose conditional mean is an unknown constant and whose conditional $(1+\alpha)$-th absolute central moment is uniformly bounded. Thus, the observations may be dependent and their conditional distributions may vary over time. The traditional i.i.d.\ model is a special case; see Section~\ref{problem set-up}.

$\Diamond$ \ We use the sharper influence-function calibration proposed in \cite{NOC} in the conditional Catoni-type framework. The resulting $\alpha$-Catoni-type supermartingales and underlying non-stitched CS form the basis of our width analysis for distributions whose conditional central $(1+\alpha)$-moment is uniformly bounded; see Section~\ref{Catoni influence function} and Lemma~\ref{catoni supermartingales}.

$\Diamond$ \ A non-asymptotic width analysis is developed for this non-stitched CS. Using the sharper constant $C_\alpha$ and the auxiliary tuning parameters $\lambda$ and $u$, we derive an explicit fixed-time upper bound on its width. At the benchmark choice $u=u_0:=4^{-1/\alpha}$ and $\lambda=\lambda_0:=\alpha$, this bound is strictly sharper than the corresponding bound in \cite[Theorem 13]{CS} for $0<\alpha<1$, whenever both sufficient conditions hold. The analysis also clarifies how the deterministic sequence $\{\theta_t\}$ and the auxiliary parameters jointly govern the effectiveness and sharpness of the bound; see Theorem~\ref{width of CS} and Remark~\ref{choose of omega}. A finite-variance comparison is also given in Corollary~\ref{variance CS} and Remark~\ref{rem:finite-variance-comparison}.

\begingroup
\setlength{\abovedisplayskip}{6pt}
\setlength{\belowdisplayskip}{6pt}
\setlength{\abovedisplayshortskip}{6pt}
\setlength{\belowdisplayshortskip}{6pt}
$\Diamond$ \ At the stitched level, applying the standard stitching method of \cite{HR2021} to the non-stitched construction yields a CS with a high-probability time-uniform width bound beyond a computable deterministic threshold, with shrinkage rate
\[
    \mathcal O\left(\left(\frac{\log\log t}{t}\right)^{\frac{\alpha}{1+\alpha}}\right)
    \quad\text{as } \ t\to\infty \ \text{ for fixed } \ \delta,
\]
and the deterministic bound expression has growth rate
\[
    \mathcal O\left((\log(1/\delta))^{\frac{\alpha}{1+\alpha}}\right)
    \quad\text{as } \ \delta\to0 \ \text{ for fixed } \ t.
\]
The width guarantee applies over the range in which the sufficient condition in Theorem~\ref{stitching CS} holds.
\par\endgroup

For each fixed $\delta$, the width of the same CS also converges to zero almost surely at the stated iterated-logarithmic rate; see Remark~\ref{rem:width-conditions}. At the benchmark choice $(u,\lambda)=(u_0,\lambda_0)$ and for $0<\alpha<1$, comparison with the Wang--Ramdas time-uniform width bound proved here shows that our bound has a strictly smaller width constant without increasing the sufficient sample-size threshold; see Remark~\ref{choose of omega in SM}. Under the bounded $(1+\alpha)$-moment assumption, this construction yields an iterated-logarithmic analogue of the benchmark fixed-time rate $\left(\log(1/\delta)/t\right)^{\alpha/(1+\alpha)}$ for mean estimation, up to constants depending on the moment bound; see \cite[Theorem~3.1]{DL2016} and \cite{LM2019}.

$\Diamond$ \ We further examine the joint asymptotic behavior of the high-probability time-uniform width bound along deterministic sequences $t_n\to\infty$ and $\delta_n\to0$, with the other parameters fixed. This bound vanishes precisely when $\log(1/\delta_n)=o(t_n)$, which also gives convergence in probability of the corresponding widths. If, in addition, $\Sigma_{n=1}^{\infty}\,\delta_n<\infty$, these widths converge to zero almost surely along the same deterministic sequences on the same data stream. The analysis also identifies two transition scales. The first, $\log(1/\delta_n)=\widetilde{\mathcal O}(\log\log t_n)$, marks the balance between the iterated-logarithmic and confidence-level terms. At the second, $\log(1/\delta_n)=\widetilde{\mathcal O}(t_n)$, the bound is of constant order, provided that the explicit sufficient condition for the width bound holds eventually. If $\log(1/\delta_n)/t_n\to\infty$, this condition eventually fails; see Remark~\ref{joint dependence}.

$\Diamond$ \ The theoretical analysis is complemented by simulation results for two heavy-tailed models with infinite variance, namely, a centered Pareto distribution with parameter $1.8$ and a Student $t$-distribution with $2$ degrees of freedom. A natural comparison is between the shrinkage and growth behavior of our improved $\alpha$-Catoni-type CS and its stitched version and those of the $\alpha$-Catoni-type CS in \cite{CS}, using their respective deterministic calibrations under common model parameters. In the reported experiments, the procedures proposed here yield smaller empirical mean interval widths. Details are given in Section~\ref{experiments}.

%% file: section/Problemset-up.tex
\section{Problem set-up and preliminaries}\label{problem set-up and preknowledge}

\subsection{Problem set-up}\label{problem set-up}~\\

We now turn to a more general setting. Let $(\Omega,\mathcal{F},(\mathcal{F}_t)_{t=1}^{\infty},\mathsf{P})$ be a filtered probability space, and let $X := (X_t)_{t=1}^{\infty}$ be a real-valued stochastic process adapted to the filtration $\{\mathcal{F}_t\}_{t=1}^{\infty}$. As usual, we take $\mathcal{F}_0$ to be the trivial $\sigma$-algebra,  $\mathcal{F}_0=\sigma(\varnothing)$.

All relations stated below, whether equalities or inequalities, that involve conditional expectations or random quantities, starting with \eqref{mean assumption}, \eqref{moment assumption}, and those that follow, are assumed to hold almost surely. Our analysis is based on the following two core assumptions, denoted by (A1) and (A2):

\textbf{(A1)} The process $X$ has a constant, unknown conditional expected value:
\begin{equation}\label{mean assumption}
  \quad \mathsf{E}[X_{t}\,|\, \mathcal{F}_{t-1}]=\mu \ \mbox{ for all } \ t \in \mathbb{N}.
\end{equation}

\textbf{(A2)} The process $X$ is conditionally  $(1 + \alpha)$-integrable for some $\alpha \in (0,1]$ and maintains a uniformly upper-bounded conditional $(1 + \alpha)$-th absolute central moment by a known constant, a real number $\nu _\alpha > 0$:
\begin{equation}\label{moment assumption}
  \mathsf{E}[|X_{t}-\mu|^{1 + \alpha}\,|\,\mathcal{F}_{t-1}] \leqslant \nu _\alpha \ \mbox{ for all } \ t \in \mathbb{N}.
\end{equation}

Our task in this paper is to construct a confidence sequence  $\{\mathrm{CI}_{t}\}_{t=1}^{\infty } $ for  $\mu$  based on the observations  $X_{1}, X_{2}, \ldots$. Recall that, for any $t \in \mathbb{N}$, \ ${\rm CI}_t = {\rm CI}_t(\omega)$ is a random interval whose endpoints depend on $X_1, \ldots, X_t$. Thus the natural requirement, written in its detailed form, is:
\begin{equation*}
  \mathsf{P}\big[\,\omega: \ \mu \in \mathrm{CI}_{t}(\omega) \ \mbox{ for all } \ t \in \mathbb{N}\,\big] \ \geqslant \ 1 - \delta.
\end{equation*}

This setting includes the traditional i.i.d.\ case mentioned in Section \ref{CS} as a special case, while also allowing one to deal with dependent observations under the same conditional-mean framework used in \cite{CS}, thereby covering a wider range of scenarios.

\begin{example}\label{example}
  Let $0<\alpha<1$, $c>0$, and let $\{\varepsilon_t\}_{t=1}^{\infty}$ be an i.i.d. sequence of Student $t$ random variables with two degrees of freedom. Define
  \begin{equation*}
    X_t = \begin{cases}
      \mu+c\,\varepsilon_1,                                           & t=1,         \\
      \mu+c\bigl(1+\mathbf{1}_{\{X_{t-1}>\mu\}}\bigr)\,\varepsilon_t, & t\geqslant2,
    \end{cases}
  \end{equation*}
  where $\mathbf{1}_A(\cdot)$ denotes the indicator of $A$. Take the natural filtration $\mathcal{F}_t=\sigma(X_1,\ldots,X_t)$ for $t\geqslant1$, with $\mathcal{F}_0$ trivial. The scale multiplying $\varepsilon_t$ is $\mathcal{F}_{t-1}$-measurable, and $\varepsilon_t$ is independent of the past and has zero mean. Thus $\mathsf{E}[X_t\mid\mathcal{F}_{t-1}]=\mu$, verifying (A1).

  Since
\begin{equation*}
  \mathsf{E}|\varepsilon_t|^{1+\alpha}
  =2^{\frac{1+\alpha}{2}} \ \frac{\Gamma((2+\alpha)/2)\Gamma((1-\alpha)/2)}{\sqrt{\pi}},
\end{equation*}
and the scale is at most $2c$, assumption (A2) holds with
\begin{equation*}
  \nu_\alpha=(2c)^{1+\alpha}\,2^{\frac{1+\alpha}{2}}\,
  \frac{\Gamma((2+\alpha)/2)\Gamma((1-\alpha)/2)}{\sqrt{\pi}}.
\end{equation*}

  Each $X_t$ nevertheless has infinite variance, since $\mathsf{E}\varepsilon_t^2=\infty$ and the predictable scale is bounded below by $c>0$. Moreover, for $t\geqslant2$, the conditional scale of $X_t-\mu$ is equal to $2c$ on the set $\{X_{t-1}>\mu\}$ and to $c$ on its complement. Both events have positive probability and give distinct conditional distributions, so the observations are dependent. This example illustrates a heavy-tailed dependent process $\{X_t\}$ for which the moment bound $\nu_\alpha$ can be specified explicitly from the model parameters.
\end{example}

The above assumptions are closely tied to the standard martingale analysis framework discussed by \cite{SY2012} and are consistent with the conditional setup considered in \cite{CS}.
Specifically, (A1) essentially requires that the sequence $\{X_t-\mu\}_{t=1}^{\infty}$ forms a martingale difference sequence, ensuring the cumulative sum $\{\Sigma_{k = 1}^{t}\, (X_k-\mu)\}_{t=1}^{\infty }$ constitutes a martingale. This structure commonly arises in scenarios involving non-i.i.d. and state-dependent noise, which frequently occur in optimization, control theory, and financial modeling contexts, as demonstrated by \cite{KY1997}.

We emphasize that the upper bound $\nu_\alpha$ in (A2) is taken throughout the paper as a {\bf{\emph{prior}}} input. This is not merely a technical convenience: in the present Catoni-type construction, $\nu_\alpha$ appears directly in the supermartingale compensator term; see Lemma~\ref{catoni supermartingales} and the discussion thereafter. More generally, some form of prior moment control is typically indispensable for uniformly valid two-sided non-asymptotic inference over broad heavy-tailed classes. In particular, in light of the Bahadur--Savage-type impossibility phenomenon \citep[Theorem~1]{BS1956}, one cannot expect nontrivial uniformly valid two-sided inference for the mean over such classes in the absence of prior moment control. Furthermore, fixed-time lower bounds over classes indexed by a $(1+\alpha)$-moment bound show that the benchmark estimation error, and hence the attainable interval width, must deteriorate as that bound increases; see \cite[Theorem~3.1]{DL2016}. In model-based settings, $\nu_\alpha$ may be derived from structural parameters, as illustrated in Example~\ref{example}. Otherwise, we simply assume that a valid upper bound is available.

\subsection{Influence function}\label{Catoni influence function}~\\

The influence function introduced by \cite{Catoni} is a crucial tool in our analysis. It is especially useful in estimating the CS width.
The high-probability width bounds established below have confidence-level dependence of order $\mathcal{O}\big((\log (1 / \delta))^{\alpha/(1+\alpha)}\big)$ over the range of sample sizes $t$ and error levels $\delta$ for which their sufficient conditions hold. This is a notably tighter dependence than the $\mathcal{O}\big((1/ \delta)^{1/(1+\alpha)}\big)$ result obtained using only Markov and Jensen inequalities. For clarity, let us provide more details.

Given a constant $C_{\alpha} > 0$ that depends only on $\alpha$, let $\psi_{\alpha} : \mathbb{R} \mapsto \mathbb{R}$ be a nondecreasing function satisfying
the following relation: \ for all $x \in \mathbb{R}$,
\begin{equation}\label{influence function}
  -\log \big(1-x+C_{\alpha}|x|^{1+\alpha}\big) \ \leqslant \ \psi_{\alpha}(x) \ \leqslant \ \log \big(1+x+C_{\alpha}|x|^{1+\alpha}\big).
\end{equation}

A function $\psi_{\alpha}$, if it exists, is called an $\alpha $-{\bf\emph{Catoni-type influence function}}.
As in \cite{Catoni}, for $\theta>0$, consider the estimating equation
\begin{equation*}
  \sum_{i=1}^{t} \psi _\alpha \big(\theta(X_{i}-\widehat{\mu}_{ \alpha})  \big) = 0, \quad t \in \mathbb{N}.
\end{equation*}

The left-hand side defines the function \  $m\mapsto\Sigma_{i=1}^{t}\,\psi_\alpha\bigl(\theta(X_i-m)\bigr)$ of the candidate mean $m$. Since $\theta>0$, increasing $m$ decreases each argument $\theta(X_i-m)$. Since $\psi_\alpha$ is nondecreasing, this function is therefore nonincreasing in $m$. Whenever the estimating equation admits a solution, any such solution defines an $\alpha$-Catoni {\rm M}-estimator $\widehat{\mu}_\alpha$. Neither existence of an exact zero nor uniqueness is required for the CS construction below.

\begingroup\setlength{\emergencystretch}{1em}
Notice that the traditional empirical mean, which corresponds to the case $\psi_\alpha(x)=x$, is significantly influenced by large $x$ values, and may fail to meet the CI width bound of order $\mathcal{O}\big((\log (1/\delta))^{\alpha/(1+\alpha)}\big)$.
This may happen when dealing with heavy-tailed distributions. \cite{LM2019} provided examples of distributions that yield a growth rate of  $\mathcal{O}\big(\sqrt{1/ \delta }\big)$ for the CI width if $\alpha = 1$, based on an i.i.d. assumption where the underlying distribution has a finite variance but lacks sub-Gaussian tail behavior.
Our analysis is carried out in the conditional martingale-difference framework and does not require independence, while yielding the logarithmic confidence-level dependence described above whenever the corresponding sufficient conditions hold.
The bounds in (\ref{influence function}) make the influence function  approximately linear near zero and restrict its magnitude to at most logarithmic growth in the tails, for large $|x|$, thus reducing the impact of extreme observations.
\par\endgroup

\begin{lemma}[\cite{NOC}]\label{C alpha}
  A nondecreasing function $\psi_\alpha$ satisfying {\rm (\ref{influence function})} exists if and only if the following relation is satisfied:
  \begin{equation}\label{c alpha}
    C_{\alpha} \ \geqslant \  \bigg(\frac{\alpha}{1+\alpha}\bigg)^{\frac{1+\alpha}{2}}\bigg(\frac{1-\alpha}{\alpha}\bigg)^{\frac{1-\alpha}{2}}.
  \end{equation}
\end{lemma}

Lemma \ref{C alpha} provides a necessary and sufficient condition for selecting the coefficient in the influence function. It is important to note that there are numerous viable choices for $C_{\alpha}$ that satisfy condition (\ref{c alpha}).
For instance,  \cite{Chen}, inspired by Taylor-like expansions, opted for $ C_{\alpha}=\frac{1}{1+\alpha} $, in which case there is a nondecreasing function, say  $\psi_{\alpha}^{\mathrm{(C)}}$,
such that
\begin{equation}\label{influence function of Chen}
  -\log \big(1-x+\frac{|x|^{1+\alpha }}{1+\alpha}\big) \ \leqslant \ \psi _\alpha^{\mathrm{(C)}} (x) \ \leqslant \ \log \big(1+x+\frac{|x|^{1+\alpha }}{1+\alpha}\big) .
\end{equation}

On the other hand,  \cite{Minsker} proposed  $C_{\alpha}= \max\Big\{\frac{\alpha}{1+\alpha}, \ \sqrt{\frac{1-\alpha}{1+\alpha}} \Big\}$,
which is a more efficient coefficient than the choice of \cite{Chen}, but is still not minimal among the admissible choices.
In our study, we adopt $C_\alpha=(\frac{\alpha}{1+\alpha})^{(1+\alpha)/2}(\frac{1-\alpha}{\alpha})^{(1-\alpha)/2}$, which is the smallest admissible value allowed by Lemma~\ref{C alpha}. With continuous extension at $\alpha=1$, this recovers the coefficient used by \cite{Catoni}, namely $C_1=1/2$. This calibration underlies the supermartingales in Lemma~\ref{catoni supermartingales} and the subsequent width analysis.

\begin{figure}[htbp]
  \centering
  \captionsetup{size=footnotesize, skip=4pt}
  \includegraphics[width=0.6\textwidth]{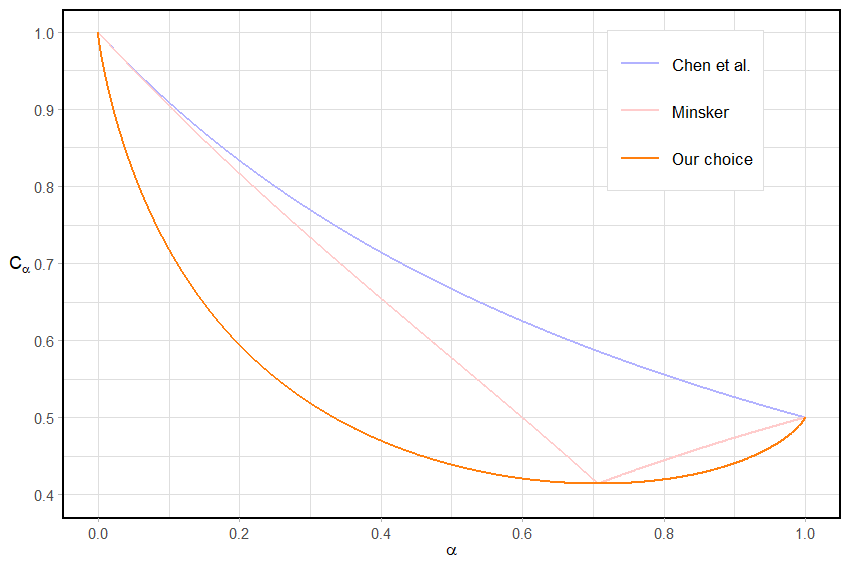}
  \caption{Comparison of different choices of \( C_{\alpha} \) as a function of \( \alpha \).}
  \label{fig:C_alpha}
\end{figure}

Figure~\ref{fig:C_alpha} illustrates how different selections of \( C_{\alpha} \) vary with the parameter \( \alpha \), comparing the choices of \cite{Chen} and \cite{Minsker} with the minimum admissible coefficient used here. Throughout, $\psi_\alpha$ and the comparison influence function $\psi_\alpha^{\mathrm{(C)}}$ are assumed to be nondecreasing and to satisfy their respective bounds~(\ref{influence function}) and~(\ref{influence function of Chen}).

\FloatBarrier

\FloatBarrier

\subsection{Nonnegative supermartingales}\label{nonnegative supermartingales} ~\\

The use of nonnegative supermartingales is a standard ingredient in modern CS constructions. Such an idea was highlighted by  \cite{HR2020}. In this context, we recall briefly the concept of a supermartingale and present two key lemmas on it.

An integrable stochastic process $\{M_t\}_{t=0}^{\infty}$, adapted to the filtration $\{\mathcal{F}_t\}_{t=0}^{\infty}$, is called a {\bf \emph{supermartingale}} if $\mathsf{E}[M_t \mid \mathcal{F}_{t-1}] \leqslant M_{t-1}$ for any $t \in \mathbb{N}$. In our applications below, we will take $M_0=1$. Since some of the supermartingales we aim to construct and use are of an exponential, multiplicative form, we need to frequently rely on the following lemma.

\begin{lemma}\label{exponential supermartingale} 
  Let $M_0=1$ and for $t\geqslant 1$, $M_{t}=\prod_{k=1}^{t} m_{k}$, where each $m_k$ is an integrable, $\mathcal{F}_k$-measurable random variable with $m_k>0$ almost surely. Then, the process $\{M_t,\mathcal{F}_t\}_{t=0}^{\infty}$ is a nonnegative supermartingale if and only if \  $\mathsf{E}[m_t\mid\mathcal{F}_{t-1}]\leqslant1$ for every $t\in\mathbb{N}$.
\end{lemma}

An essential feature of the nonnegative supermartingales is the validity of what is known as Ville's inequality, first established in \cite{Ville}. It  is based on an idea behind Markov's inequality and expands it to apply over an unlimited time period.

\begin{lemma}[Ville's Inequality]\label{Ville's Inequality} 
  Let  $\{M_{t}, \mathcal{F}_t\}_{t=0}^{\infty } $ be a nonnegative supermartingale with $ M_{0}=1 $. Then for any  $\delta  \in(0,1)$, we have the
  following:
  \begin{equation*}
    \mathsf{P}\big[\,\omega: \ \exists \ t \in \mathbb{N} \ \text{ such that } \ M_{t}(\omega) \geqslant 1 / \delta \big] \ \leqslant \ \delta .
  \end{equation*}
\end{lemma}

A concise proof of Lemma \ref{Ville's Inequality} is available in  \cite{HR2021}.

Ville's inequality provides the standard supermartingale-to-CS conversion that underlies our construction. This standard device, together with the properties of nonnegative supermartingales, forms the backbone of the sequential framework used below.
The following supermartingales have the same form as those in \cite{FangCS}, when the influence functions and tuning sequences are matched. We give the verification under (A1)--(A2).

\begin{lemma}[Improved $\alpha$-Catoni supermartingales]\label{catoni supermartingales}
  Assume (A1) and (A2). Let $\{\theta_t,\mathcal{F}_t\}_{t=0}^{\infty}$ be any predictable process such that
  $\theta_t \geqslant 0$ for all $t \geqslant 0$, and let $\psi_\alpha$ be an $\alpha$-Catoni-type
  influence function satisfying condition {\rm (\ref{influence function})}. Define two processes,
  $\{M_t^{+}\}_{t=0}^{\infty}$ and $\{M_t^{-}\}_{t=0}^{\infty}$, as follows:
  \begin{equation}\label{eq:centered-supermartingales}
    \begin{aligned}
    M_t^{+} & =
    \begin{cases}
      1,                                                                                                              & t = 0,            \\
      \prod\limits_{i=1}^{t} \exp \big\{ \psi_\alpha(\theta_i(X_i-\mu))-\theta_i^{1+\alpha}C_\alpha\nu_\alpha \big\}, & t \in \mathbb{N},
    \end{cases} \\[4pt]
    M_t^{-} & =
    \begin{cases}
      1,                                                                                                              & t = 0,            \\
      \prod\limits_{i=1}^{t} \exp \big\{-\psi_\alpha(\theta_i(X_i-\mu))-\theta_i^{1+\alpha}C_\alpha\nu_\alpha \big\}, & t \in \mathbb{N}.
    \end{cases}
    \end{aligned}
  \end{equation}

  Then both $\{M_t^{+}, \mathcal{F}_t\}_{t=0}^{\infty}$ and $\{M_t^{-}, \mathcal{F}_t\}_{t=0}^{\infty}$
  are nonnegative supermartingales.

  Moreover, for every fixed $m\in\mathbb{R}$ and $u\in(0,1)$, let $\widetilde M_0^\pm(m,u)=1$ and, for $t\in\mathbb{N}$, define
  \begin{equation}\label{eq:compensated-supermartingales}
    \widetilde M_t^\pm(m,u)
    =\prod_{i=1}^{t}\exp\bigg\{\pm\psi_\alpha\bigl(\theta_i(X_i-m)\bigr)
    \mp\theta_i(\mu-m)-C_\alpha\theta_i^{1+\alpha}
    \bigl[u^{-\alpha}\nu_\alpha+(1-u)^{-\alpha}|\mu-m|^{1+\alpha}\bigr]\bigg\}.
  \end{equation}

  Then $\{\widetilde M_t^\pm(m,u),\mathcal{F}_t\}_{t=0}^{\infty}$ are also nonnegative supermartingales.
\end{lemma}

\begin{proof}
  Recalling that $\mathcal{F}_0$ is the trivial $\sigma$-algebra, we note that the initial values $M_0^+ = M_0^- = 1$ are $\mathcal{F}_0$-measurable.

  For every  $t \in \mathbb{N} $, we see that the following inequality holds:
  $$\begin{aligned}
      \mathsf{E}[\,\exp \{ \psi _\alpha  (\theta _{t}(X_{t}-\mu))\} \,|\, \mathcal{F}_{t-1}] \leqslant & ~\mathsf{E}[ 1+\theta _{t}(X_{t}-\mu)+ C_{\alpha}\theta _{t}^{1+\alpha}|X_{t}-\mu|^{1+\alpha} \mid \mathcal{F}_{t-1}] \\
      \leqslant                                                                                        & ~1+C_{\alpha} \nu _\alpha \theta _{t}^{1+\alpha}  \leqslant  \exp \{C_{\alpha} \nu _\alpha \theta _{t}^{1+\alpha}\}.
    \end{aligned}$$

  This can be equivalently expressed as:
  $$ \mathsf{E}\,[\exp \{ \psi _\alpha  (\theta _{t}(X_{t}-\mu))-\theta _{t}^{1+\alpha } C_{\alpha}\nu _\alpha\} \mid \mathcal{F}_{t-1}] \ \leqslant \ 1 .$$

  Using the lower bound in (\ref{influence function}) gives the same bound with $-\psi_\alpha$ in place of $\psi_\alpha$.

  Next, fix $m\in\mathbb{R}$ and $u\in(0,1)$. By convexity of $x\mapsto x^{1+\alpha}$ on $[0,\infty)$, we have
  \begin{equation}\label{convex-splitting-u}
    (a+b)^{1+\alpha}
    \leqslant u^{-\alpha}a^{1+\alpha}+(1-u)^{-\alpha}b^{1+\alpha} \quad \mbox{ for any }
    \quad a,b\geqslant0.
  \end{equation}

  Applying this inequality with $a=|X_t-\mu|$ and $b=|\mu-m|$, and using (A1) and (A2), gives
  \[
    \begin{aligned}
      \mathsf{E}[X_t-m\mid\mathcal{F}_{t-1}] & =\mu-m,                                                            \\
      \mathsf{E}[|X_t-m|^{1+\alpha}\mid\mathcal{F}_{t-1}]
                                             & \leqslant u^{-\alpha}\nu_\alpha+(1-u)^{-\alpha}|\mu-m|^{1+\alpha}.
    \end{aligned}
  \]

  Therefore, predictability of $\theta_t$, the two bounds in (\ref{influence function}), and $1+x\leqslant\exp(x)$ yield
  \[
    \begin{aligned}
       & \mathsf{E}\bigl[\exp\{\pm\psi_\alpha(\theta_t(X_t-m))\}\mid\mathcal{F}_{t-1}\bigr] \\
       & \quad\leqslant1\pm\theta_t(\mu-m)
      +C_\alpha\theta_t^{1+\alpha}
      \bigl[u^{-\alpha}\nu_\alpha+(1-u)^{-\alpha}|\mu-m|^{1+\alpha}\bigr]                   \\
       & \quad\leqslant\exp\Big\{\pm\theta_t(\mu-m)
      +C_\alpha\theta_t^{1+\alpha}
      \bigl[u^{-\alpha}\nu_\alpha+(1-u)^{-\alpha}|\mu-m|^{1+\alpha}\bigr]\Big\}.
    \end{aligned}
  \]

  Multiplication by the predictable compensating exponential shows that each factor of $\widetilde M_t^\pm(m,u)$ has conditional expectation at most $1$. The factors of all four processes are strictly positive and integrable, so Lemma~\ref{exponential supermartingale} proves the claims.
\end{proof}

We next clarify the role of condition (A2) in Lemma \ref{catoni supermartingales}. In the present construction, the conditional $(1+\alpha)$-moment bound is used in an essential way to verify the supermartingale property of the processes $\{M_t^+, \mathcal F_t\}_{t=0}^\infty$ and $\{M_t^-, \mathcal F_t\}_{t=0}^\infty$; see \cite[Proposition 3]{CS}.
In this sense, the compensation term of the centered processes $M_t^\pm$ in Lemma \ref{catoni supermartingales} is calibrated as {\bf \emph{sharply}} as the argument allows within the current Catoni-type framework. The compensated processes $\widetilde M_t^\pm(m,u)$ will be used in the width analysis below. At $m=\mu$, their moment compensator retains the additional factor $u^{-\alpha}$. If condition (A2) fails to hold on a set of positive probability, then the above supermartingale verification may break down, and the resulting time-uniform guarantee can no longer be justified by this route.

\subsection{Stitching method}\label{stitching method}~\\ 

Recall that the {\bf \emph{stitching method}} is a standard technique for combining CIs or CSs that are valid over different ranges of sample sizes; see \cite{HR2021}. In the finite-variance case, \cite{CS} used this technique to obtain a law-of-the-iterated-logarithm-type upper bound of order $\mathcal{O}\big(\sqrt{(\log \log t)/t}\big)$ for the CS width; see also \cite{DH2024} for related lower-bound results in regular models.

In the present heavy-tailed setting, the same construction leads to a stitched CS whose width admits a high-probability time-uniform upper bound of order
\[
  \mathcal{O}\big(((\log \log t)/t)^{\frac{\alpha}{1+\alpha}}\big)
\quad \text{ as } \ t\to\infty \ \text{ for fixed } \ \delta.
\]

Notice that Theorem~\ref{stitching CS} establishes this time-uniform width bound beyond an explicit deterministic threshold.

We now describe the main steps of the stitching method.

\begin{lemma}[Stitching method]\label{Stitching method} 
  Consider a partition $\{T_k\}_{k = 0}^{\infty}$ of the index set $\mathbb{N}$ (that is, $\mathbb{N}=\bigcup_{k=0}^{\infty} T_k$, where $T_i\cap T_j= \varnothing$ for all $i\neq j$), and a sequence of positive `small' numbers $\{\delta _k\}_{k = 0}^{\infty}$ such that \  $\Sigma_{k = 0}^{\infty}\,\delta _k \leqslant \delta$.
  For each $k\in\mathbb N_0$, assume the existence of a $(1-\delta_k )$-{\rm CS}, say $\mathrm{CI}_{t}^{(k)}$, over the index set $T_k$, for the unknown $\mu$:
  \[
    \mathsf{P}\big[\,\omega:\ \mu \in \mathrm{CI}_{t}^{(k)}(\omega) \ \text{ for all } \ t \in T_k\big] \geqslant 1-\delta_k.
  \]

  Then, the stitched-{\rm CS}, say ${\rm CS}^{\rm stch}=\{\mathrm{CI}_{t}^{\mathrm{stch}},\ t\in\mathbb{N}\}$, defined by
  \[
    \mathrm{CI}_{t}^{\mathrm{stch}}=\mathrm{CI}_{t}^{(k)}, \quad \text{ for } \ t \in T_k,\ k\in\mathbb N_0,
  \]
  represents a $(1-\delta )$-{\rm CS} for $\mu$.
\end{lemma}

\begin{proof}
  Let us rewrite the involved random sets as unions:
  \[
    \big\{\,\exists\, t \in \mathbb{N} \ \text{ such that } \ \mu \notin \mathrm{CI}_{t}^{\mathrm{stch}}\big\}
    =
    \bigcup_{t\in \mathbb{N}}\big\{\mu \notin \mathrm{CI}_{t}^{\mathrm{stch}}\big\}
    =
    \bigcup_{k=0}^{\infty}\bigcup_{t\in T_k}\big\{\mu \notin \mathrm{CI}_{t}^{(k)}\big\}.
  \]

  Given that $\mathrm{CI}_{t}^{(k)}$ is a $(1-\delta_k )$-{\rm CS} in $T_k$ for $\mu$, we have
  \[
    \mathsf{P}\big[\,\exists\, t \in T_k \ \text{ such that } \ \mu \notin \mathrm{CI}_{t}^{(k)}\big] \leqslant \delta_k.
  \]

  Combining the above, we deduce
  \[
    \begin{aligned}
                & ~\mathsf{P}\big[\,\exists\, t \in \mathbb{N} \ \text{ such that }\ \mu \notin \mathrm{CI}_{t}^{\mathrm{stch}}\big]
      = \mathsf{P}\left[\bigcup_{k=0}^{\infty}\bigcup_{t\in T_k}\{\mu \notin \mathrm{CI}_{t}^{(k)}\}\right]                      \\
      \leqslant & ~\sum_{k=0}^{\infty}\mathsf{P}\left[\bigcup_{t\in T_k}\{\mu \notin \mathrm{CI}_{t}^{(k)}\}\right]
      = \sum_{k=0}^{\infty}\mathsf{P}\left[\,\exists\, t \in T_k \ \text{ such that } \ \mu \notin \mathrm{CI}_{t}^{(k)}\right]       \\
      \leqslant & ~\sum_{k=0}^{\infty}\delta_k \leqslant \delta.
    \end{aligned}
  \]

  Hence, the stitched-CS, namely $\{\mathrm{CI}_{t}^{\mathrm{stch}}\}_{t=1}^{\infty}$, is indeed a $(1-\delta )$-CS for the unknown $\mu$.
\end{proof}

In comparison with the original CS construction in Definition \ref{CS definition}, the stitched-CS provides a refined way to reduce the width through an appropriate choice of both the partition $\{T_k\}_{k=0}^{\infty}$ and the sequence $\{\delta_k\}_{k=0}^{\infty}$.

A finer partition allows the tuning coefficients to be tailored to shorter time ranges, while the error probabilities assigned to these ranges must still sum to at most $\delta$. The potential gain from more localized tuning must be balanced against the larger terms involving the error levels in the width bounds associated with smaller error allocations. Hence, the choice of both $\{T_k\}_{k=0}^{\infty}$ and $\{\delta_k\}_{k=0}^{\infty}$ is crucial. A concrete specification of these sequences is given in Section~\ref{main results}, Theorem~\ref{stitching CS}.

%% file: section/Mainresults.tex
\section{Main results and proofs}\label{main results}
Combining the two centered supermartingales in Lemma~\ref{catoni supermartingales} with the standard supermartingale-to-CS conversion given by Ville's inequality (Lemma~\ref{Ville's Inequality}), we obtain our first main result in this section, namely, an improved $\alpha$-Catoni-type CS for $\mu$. Since this argument is also used in the subsequent width analysis, we provide a short proof for completeness.

\begin{theorem}[Improved $\alpha$-Catoni-type CS]\label{Catoni-type confidence sequence} 
  Assume that conditions {\rm (A1)} and {\rm (A2)} hold. Let $\{\theta_t\}_{t=1}^{\infty}$ be any nonnegative predictable process, and $\psi_\alpha$ a nondecreasing $\alpha$-Catoni-type influence function satisfying condition {\rm (\ref{influence function})}.

  The following random intervals $\{\mathrm{CI}_t\}_{t=1}^{\infty}$ constitute a $(1-\delta)$-{\rm CS} for the unknown $\mu$:
  \begin{equation}\label{Catoni-type CS}
    \mathrm{CI}_{t}
    =
    \bigg\{
    m \in \mathbb{R}: \
    \Big\lvert
    \sum_{i=1}^{t} \psi_\alpha\big(\theta_i(X_i-m)\big)
    \Big\rvert
    \leqslant
    C_{\alpha}\nu_{\alpha}\sum_{i=1}^{t}\theta_i^{1+\alpha}
    +
    \log\Big(\frac{2}{\delta}\Big)
    \bigg\}.
  \end{equation}
\end{theorem}

\begin{proof}
  By Lemma \ref{catoni supermartingales}, the centered processes $M_t^+$ and $M_t^-$ in \eqref{eq:centered-supermartingales} are nonnegative supermartingales with initial value $1$. Applying Ville's inequality (Lemma \ref{Ville's Inequality}) to these two processes yields
  \[
    \mathsf{P}\bigg(
    \exists ~ t \in \mathbb{N}: \
    \sum_{i=1}^t \psi_\alpha\big(\theta_i(X_i-\mu)\big)
    -
    C_\alpha \nu_\alpha \sum_{i=1}^t \theta_i^{1+\alpha}
    \geqslant \log\Big(\frac{2}{\delta}\Big)
    \bigg)
    \leqslant \frac{\delta}{2},
  \]
  and
  \[
    \mathsf{P}\bigg(
    \exists ~ t \in \mathbb{N}: \
    -\sum_{i=1}^t \psi_\alpha\big(\theta_i(X_i-\mu)\big)
    -
    C_\alpha \nu_\alpha \sum_{i=1}^t \theta_i^{1+\alpha}
    \geqslant \log\Big(\frac{2}{\delta}\Big)
    \bigg)
    \leqslant \frac{\delta}{2}.
  \]

  Therefore, by the union bound, with probability at least $1-\delta$, for all $t \in \mathbb{N}$,
  \[
    \Big\lvert
    \sum_{i=1}^t \psi_\alpha\big(\theta_i(X_i-\mu)\big)
    \Big\rvert
    \leqslant
    C_\alpha \nu_\alpha \sum_{i=1}^t \theta_i^{1+\alpha}
    +
    \log\Big(\frac{2}{\delta}\Big).
  \]

  Equivalently, with probability at least $1-\delta$, one has $\mu\in \mathrm{CI}_t$ for all $t\geqslant 1$. Hence, $\{\mathrm{CI}_t\}_{t=1}^\infty$ is indeed a $(1-\delta)$-{\rm CS} for $\mu$.
\end{proof}

\begin{remark}[Predictable tuning]
  In \eqref{Catoni-type CS}, the coefficient $\theta_t$ may depend on past observations, provided that it is nonnegative and $\mathcal F_{t-1}$-measurable. A more detailed interpretation of the deterministic tuning sequence $\{\theta_t\}_{t=1}^{\infty}$, together with the roles of the auxiliary parameters $C_\alpha$, $u$, and $\lambda$ in the width analysis, is given in Remark~\ref{choose of omega}.
\end{remark}

Despite the fact that this CS does not admit a closed-form expression, its endpoints can be determined by inverting the inequalities in \eqref{Catoni-type CS} using monotone threshold searches, since the function $m \mapsto \Sigma_{i=1}^{t}\, \psi_\alpha\big(\theta_i(X_i-m)\big)$ is nonincreasing. We next derive an upper bound for the width of the improved $\alpha$-Catoni-type CS.

In what follows, we use the following notation:
\begin{equation}\label{eq:log-error-budget}
  L(\beta, \delta) := \log (2/\beta) + \log (2/\delta).
\end{equation}

\begin{theorem}\label{width of CS} 
  Assume that conditions {\rm (A1)} and {\rm (A2)} hold, and that the coefficients $\{\theta _{t}\}_{t=1}^{\infty }$ are nonrandom and nonnegative. Consider a time $t$ such that $\Sigma_{i=1}^t\,\theta_i>0$. Let $\beta\in(0,1)$, $0<\lambda\leqslant 1/\alpha$, and $0<u<1$ be chosen such that
  \begin{equation}\label{condition}
    C_{\alpha }^{-\frac{1}{\alpha } }(1-u )\frac{\lambda^{\frac{1}{\alpha}}}{(1+\lambda)^{1+\frac{1}{\alpha }}}\bigg(\frac{\sum_{i=1}^{t} \theta _i}{\sum_{i=1}^{t} \theta _{i}^{1+\alpha }}\bigg)^{1+\frac{1}{\alpha }} \geqslant (1+u  ^{-\alpha })C_{\alpha} \nu _\alpha+\frac{L(\beta, \delta)}{\sum_{i=1}^{t} \theta _{i}^{1+\alpha }}.
  \end{equation}

  Then, with probability at least $1-\beta$, the width  $|\mathrm{CI}_{t}|$ can be bounded as follows:
  \begin{equation}\label{result}
    |\mathrm{CI}_{t}| \leqslant  \frac{2(1+\lambda)}{\sum_{i=1}^{t} \theta_i} \bigg( (1+u^{-\alpha })C_{\alpha} \nu_\alpha\sum_{i=1}^{t} \theta _{i}^{1+\alpha }+ L(\beta, \delta)\bigg).
  \end{equation}

  Here $C_{\alpha} =\big(\frac{\alpha}{1+\alpha}\big)^{(1+\alpha)/2}\big(\frac{1-\alpha}{\alpha}\big)^{(1-\alpha)/2}$ is the constant mentioned in Section \rm{\ref{Catoni influence function}.}
\end{theorem}

The fixed-time sufficient condition \eqref{condition} and width bound \eqref{result} are a specialization of the analysis in \cite{FangCS}. We provide below a self-contained proof of Theorem~\ref{width of CS} under (A1)--(A2), which will be useful for the subsequent width comparisons and stitching analysis.

\begin{remark}[Parameter choices and width comparison]\label{choose of omega}
  First, $C_\alpha$ is not a free tuning parameter, but the structural constant associated with the influence function introduced in Section \ref{Catoni influence function}. The choice of the smallest admissible value of $C_\alpha$ yields a sharper calibration in the resulting supermartingale construction and in the corresponding width bound.

  Second, the parameters $u$ and $\lambda$ should not be viewed as playing completely parallel roles. On one hand, smaller values of $\lambda$ lead to a sharper explicit upper bound in \eqref{result}, but make the sufficient condition \eqref{condition} harder to satisfy. Values of $u$ closer to $1$ also reduce the displayed width coefficient; their effect on the sufficient condition is governed by \eqref{condition}. On the other hand, the parameter $u$ also enters the concrete deterministic choice of $\{\theta_t\}_{t=1}^{\infty}$ used later; see, in particular, \eqref{theta of me}. Thus, although both $u$ and $\lambda$ contribute to the trade-off underlying the theorem-level width bound, the parameter $u$ has a more direct influence on the subsequent implementation through the chosen sequence $\{\theta_t\}_{t=1}^{\infty}$.

  Third, the deterministic sequence $\{\theta_t\}_{t=1}^{\infty}$ determines the CS shrinkage rate. Roughly speaking, $\Sigma_{i=1}^t\, \theta_i$ drives the shrinkage of the interval, whereas $\Sigma_{i=1}^t\, \theta_i^{1+\alpha}$ reflects the heavy-tail penalty. More precisely, condition (\ref{condition}) requires a suitable combination of these two quantities to be large enough, while the upper bound in (\ref{result}) is also governed by the same two sums.
  Hence, the choice of $\{\theta_t\}_{t=1}^{\infty}$ defines these two sums which jointly determine and control both the improved bound and the resulting width.

  In particular, as shown in Table \ref{table}, when $\theta_t \asymp t^{-1/(1+\alpha)}$, one has
  $\Sigma_{i=1}^t\, \theta_i \asymp t^{\alpha/(1+\alpha)}$ and
  $\Sigma_{i=1}^t\, \theta_i^{1+\alpha} \asymp \log t$.
  Hence, as $t\to\infty$,
  \[
    \left(\frac{\sum_{i=1}^t \theta_i}{\sum_{i=1}^t \theta_i^{1+\alpha}}\right)^{1+\frac{1}{\alpha}} \
    \asymp \
    \frac{t}{(\log t)^{1+\frac{1}{\alpha}}} \ \to \ \infty
    \qquad
    \text{and}
    \qquad
    \frac{L(\beta,\delta)}{\sum_{i=1}^t \theta_i^{1+\alpha}} \
    \asymp \
    (\log t)^{-1} \ \to \ 0.
  \]

  Therefore, for this canonical deterministic choice, condition (\ref{condition}) is eventually satisfied. Similar conclusions also apply to the other deterministic choices of $\{\theta_t\}$ listed in Table \ref{table}.

  Comparing \eqref{result} with the width bound in \cite[Theorem 13]{CS} at a common deterministic sequence $\{\theta_t\}_{t=1}^{\infty}$, under the sufficient conditions of both results, shows that the benchmark choice $(u,\lambda)=(u_0,\lambda_0)$ yields a strictly sharper explicit upper bound for $0<\alpha<1$. The two bounds coincide at $\alpha=1$ under this choice. More generally, choices of $\lambda$ and $u$ can further sharpen the explicit upper bound in \eqref{result}, subject to its sufficient condition, while the sequence $\{\theta_t\}_{t=1}^{\infty}$ determines the overall shrinkage rate.
\end{remark}

\begin{proof}[Proof of Theorem \rm{\ref{width of CS}}]
  For the candidate mean $m\in\mathbb{R}$, define
  \begin{equation}\label{eq:width-score}
    g_t(m)=\sum_{i=1}^t\psi_\alpha\big(\theta_i(X_i-m)\big).
  \end{equation}

  Since $\psi_\alpha$ is nondecreasing and $\theta_i\geqslant0$, the function $g_t$ is nonincreasing in $m$. For each fixed $m\in\mathbb R$, Lemma~\ref{catoni supermartingales} implies that the processes defined in \eqref{eq:compensated-supermartingales} satisfy
  \[
    \mathsf E\big[\widetilde M_t^\pm(m,u)\big]\leqslant1.
  \]

  For convenience, denote
  \[
    U :=(1-u)^{-\alpha}|\mu-m|^{1+\alpha}.
  \]

  Because the tuning coefficients are nonrandom, the compensating terms in these processes are deterministic. We therefore obtain the following two exponential-moment bounds:
  \[
    \begin{aligned}
      \mathsf E\big[\exp(g_t(m))\big]
       & \leqslant\exp\bigg\{\sum_{i=1}^t\theta_i(\mu-m)
      +C_\alpha\sum_{i=1}^t\theta_i^{1+\alpha}(u^{-\alpha}\nu_\alpha+U)\bigg\}, \\
      \mathsf E\big[\exp(-g_t(m))\big]
       & \leqslant\exp\bigg\{-\sum_{i=1}^t\theta_i(\mu-m)
      +C_\alpha\sum_{i=1}^t\theta_i^{1+\alpha}(u^{-\alpha}\nu_\alpha+U)\bigg\}.
    \end{aligned}
  \]

  Define the functions $B_{t}^{+}(m)$ and $B_{t}^{-}(m)$ as follows:
  \begin{align}
    B_{t}^{+}(m) & =\sum_{i=1}^{t}(\theta _{i}(\mu-m)+C_{\alpha}\theta _{i}^{1+\alpha }(u ^{-\alpha }\nu_\alpha+ U
    ))+\log (2/\beta) \notag                                                                                                                                \\*
                 & =\sum_{i=1}^{t}C_{\alpha}\theta _{i}^{1+\alpha } U
    -\sum_{i=1}^{t} \theta _{i}(m-\mu)+\sum_{i=1}^{t} C_{\alpha}\theta _{i}^{1+\alpha }u ^{-\alpha }\nu _\alpha+\log (2 / \beta ) ; \notag \displaybreak[0] \\
    B_{t}^{-}(m) & =\sum_{i=1}^{t}\theta_i(\mu-m)-C_\alpha\sum_{i=1}^{t}\theta_i^{1+\alpha}(u^{-\alpha}\nu_\alpha+U)-\log (2/\beta) \notag                  \\*
                 & =-\sum_{i=1}^{t}C_{\alpha}\theta _{i}^{1+\alpha } U
    -\sum_{i=1}^{t} \theta _{i}(m-\mu)-\sum_{i=1}^{t} C_{\alpha}\theta _{i}^{1+\alpha }u ^{-\alpha }\nu _\alpha-\log (2 / \beta ) .\label{eq:width-envelopes}
  \end{align}

  By applying Markov's inequality, we find that for each fixed $m \in \mathbb{R},$
  \begin{equation}\label{upper bound inequality}
    \mathsf{P}[g_{t}(m) \geqslant B_{t}^{+}(m)] \leqslant \frac{\mathsf{E}\big[ \exp (g_{t}(m))\big]}{\exp(B_{t}^{+}(m))} \leqslant \frac{\beta}{2},
  \end{equation}
  \begin{equation*}
    \mathsf{P}[g_{t}(m) \leqslant B_{t}^{-}(m)] \leqslant \frac{\mathsf{E}\big[\exp (-g_{t}(m))\big]}{\exp(-B_{t}^{-}(m))} \leqslant \frac{\beta}{2} .
  \end{equation*}

  Thus, for each fixed $m\in\mathbb R$,
  \[
    \mathsf P\big[B_t^-(m)<g_t(m)<B_t^+(m)\big]\geqslant1-\beta.
  \]

  These bounds apply separately at fixed comparison points; no simultaneous event over all $m$ is asserted.

  With $B_t^+(m)$ as given above, consider the following equation:
  \begin{equation}\label{upper bound equation}
    B_{t}^{+}(m)=-C_{\alpha} \nu _\alpha \sum_{i=1}^{t} \theta _{i}^{1+\alpha }-\log (2 / \delta ).
  \end{equation}

  Define  $m=\rho_{t}$, where $\rho_t$ is the smallest solution (the existence of roots is detailed below) of equation (\ref{upper bound equation}). Given that the sequence $\{\theta_{t}\}$ is nonrandom, the quantity $\rho _{t}$ is also nonrandom. Substituting $m=\rho_{t}$  into (\ref{upper bound inequality}), we obtain
  $$\mathsf{P}\left[g_{t}(\rho _{t}) \geqslant-C_{\alpha} \nu _\alpha \sum_{i=1}^{t} \theta _{i}^{1+\alpha }-\log (2 / \delta )\right] \leqslant \beta  / 2 .$$

  With probability at least $1-\beta/2$, we therefore have
  \[
    g_t(\rho_t)<-C_\alpha\nu_\alpha\sum_{i=1}^t\theta_i^{1+\alpha}-\log(2/\delta).
  \]

  On this event, every $m\geqslant\rho_t$ satisfies $g_t(m)\leqslant g_t(\rho_t)$ and hence lies outside $\mathrm{CI}_t$. Consequently, using the convention $\sup\varnothing=-\infty$, we see that
  \begin{equation}\label{upper}
    \mathsf{P}[\sup\mathrm{CI}_{t} \leqslant \rho_t] \geqslant 1-\frac{\beta}{2}.
  \end{equation}

  Hence, the next step is to estimate the difference $\rho _{t}-\mu $.

  To facilitate easier calculation, let us rewrite  relation (\ref{upper bound equation}) as follows, with each  symbol $\sum$  meaning a summation  $\sum_{i=1}^{t}$:
  $$C_{\alpha}(1-u )^{-\alpha }\,\frac{\sum \,\theta _{i}^{1+\alpha }}{\sum \,\theta _{i}} \,\lvert m-\mu\rvert ^{1+\alpha } - (m-\mu) + \ C_{\alpha}\nu _\alpha(1+u ^{-\alpha })\,\frac{\sum \,\theta _{i}^{1+\alpha }}{\sum \,\theta _{i}}+\frac{L(\beta,\delta)}{\sum \, \theta _{i}}=0.$$

  Let us denote
  \begin{equation}\label{eq:width-root-coefficients}
    K_t=C_{\alpha}(1-u )^{-\alpha }\,\frac{\sum \,\theta _{i}^{1+\alpha }}{\sum \,\theta _{i}}, \quad J_t=C_{\alpha}\nu _\alpha(1+u ^{-\alpha })\,\frac{\sum \, \theta _{i}^{1+\alpha }}{\sum \, \theta _{i}}+\frac{L(\beta,\delta)}{\sum \,\theta _{i}}~,
  \end{equation}
  where $K_t>0$ and $J_t>0$. The left-hand side of the preceding equation is positive for $m\leqslant\mu$, so every solution satisfies $m>\mu$. Setting $z=m-\mu$, equation (\ref{upper bound equation}) becomes
  \begin{equation}\label{root equation I}
    K_t z^{1+\alpha }-z+J_t=0, \quad z > 0 .
  \end{equation}

  Denoting further
  $$\widetilde J_t=K_t^{\frac{1}{\alpha}} J_t,$$
  and setting  $y=K_t^{1/\alpha} z$, we  transform (\ref{root equation I}) into the equation
  \begin{equation}\label{root equation II}
    f_t(y) := y^{1+\alpha }-y+\widetilde J_t=0, \quad y > 0 .
  \end{equation}

  Note that the function  $f_t( \ \cdot \ )$ is strictly convex in $y$  and  $f_t(y) \rightarrow \infty $ as  $y \rightarrow \infty $. Therefore,  $f_t( \ \cdot \ )$ has a unique minimum on $[0,\infty)$, attained at the point $y_*$, where
  $$y_{*}=\bigg(\frac{1}{1+\alpha }\bigg)^{\frac{1}{\alpha}} \quad \Rightarrow \quad \min _{y \geqslant 0} f_t(y)= f_t (y_{*})= \widetilde J_t-\frac{\alpha}{(1+\alpha )^{1+\frac{1}{\alpha } } } .$$

  If this minimum is non-positive, we have
  \begin{equation}\label{root existence condition}
    \widetilde J_t \leqslant \frac{\alpha }{(1+\alpha )^{1+\frac{1}{\alpha } } }.
  \end{equation}

  Then the equation $ f_t(y)=0$ has one positive root, or two positive roots if the inequality (\ref{root existence condition}) is strict, because $f_t(0)>0$ and  $y_{*}>0$. Define $y_t \in [0, y_{*}]$ as the smaller of these roots. Since $y_t=K_t^{1/\alpha}(\rho _{t}-\mu )$, estimating the difference $\rho _{t}-\mu$ is tantamount to determining the value of $y_t$.

  If, instead of (\ref{root existence condition}), the following condition is satisfied for some $\lambda > 0$,
  \begin{equation}\label{parameter condition}
    \widetilde J_t \leqslant \frac{\lambda ^{\frac{1}{\alpha}}}{(1+\lambda )^{1+\frac{1}{\alpha }}} := h (\lambda ),
  \end{equation}
  then it results in the relation $ f_t((1+\lambda ) \widetilde J_t) \leqslant 0$. Note that condition (\ref{parameter condition}) aligns with condition (\ref{root existence condition}) when $\lambda =\frac{1}{\alpha}$. At this point, the function $h(\lambda)$ reaches its maximum. Thus, if condition (\ref{parameter condition}) holds, it implies that condition (\ref{root existence condition}) is also satisfied, thereby ensuring the existence of roots for equation (\ref{root equation II}).

  Thus, the positive solution $y_t$ of equation (\ref {root equation II}) satisfies the relation
  \begin{equation}\label{scaled-root-bound}
    y_t \leqslant(1+\lambda)\widetilde J_t,
    \qquad \Rightarrow \qquad
    \rho _{t}-\mu \leqslant(1+\lambda ) J_t .
  \end{equation}

  The bound in \eqref{scaled-root-bound} also holds at $\lambda=1/\alpha$, including the case $\widetilde J_t=h(1/\alpha)$, where the two roots coincide at $y_*=(1+1/\alpha)\widetilde J_t$.
  Equivalently,
  \begin{equation}\label{upper bound}
    \rho _{t}-\mu \leqslant(1+\lambda )\bigg[C_{\alpha}\nu _\alpha(1+u ^{-\alpha })\frac{\sum \, \theta _{i}^{1+\alpha }}{\sum \, \theta _{i}}+\frac{L(\beta,\delta)}{\sum \, \theta _{i}}\bigg] .
  \end{equation}

  Assume now that condition (\ref{parameter condition}) holds. It can  be reformulated as
  $$C_{\alpha}\nu_\alpha(1+u^{-\alpha })\,\frac{\sum \, \theta _{i}^{1+\alpha }}{\sum \, \theta _{i}} +  \frac{L(\beta,\delta)}{\sum \, \theta _{i}} \leqslant  \frac{\lambda ^{\frac{1}{\alpha}}}{(1+\lambda )^{1+\frac{1}{\alpha }}}C_{\alpha}^{-\frac{1}{\alpha}}(1-u)\bigg(\frac{\sum \, \theta _{i}}{\sum \, \theta _{i}^{1+\alpha }}\bigg)^{\frac{1}{\alpha}}.$$

  This coincides with  condition (\ref{condition}) after a simple rearrangement of terms.

  Combining the outcomes of (\ref{upper}) and (\ref{upper bound}), we obtain a one-sided bound:
  \begin{equation*}
    \mathsf{P}\bigg\{\sup\mathrm{CI}_{t} \leqslant \mu + (1+\lambda )\bigg[C_{\alpha}\nu _\alpha(1+u ^{-\alpha })\,\frac{\sum \, \theta _{i}^{1+\alpha }}{\sum \, \theta _{i}}+\frac{L(\beta,\delta)}{\sum \, \theta _{i}}\bigg]\bigg\} \geqslant 1-\frac{\beta}{2} .
  \end{equation*}

  Now, instead of $\rho_t$ used
  above, let $\varrho_t$ be the largest $m$ which is a solution to the equation
  \begin{equation*}
    B_{t}^{-}(m)=C_{\alpha} \nu _\alpha \sum_{i=1}^{t} \theta _{i}^{1+\alpha }+\log (2 / \delta ) .
  \end{equation*}

  This point is also nonrandom. The lower-tail Markov bound gives
  \[
    \mathsf P\bigg[g_t(\varrho_t)\leqslant C_\alpha\nu_\alpha\sum_{i=1}^t\theta_i^{1+\alpha}+\log(2/\delta)\bigg]\leqslant\frac{\beta}{2}.
  \]

  On the complementary event, $g_t(\varrho_t)$ is strictly above the upper threshold in \eqref{Catoni-type CS}. Since $g_t$ is nonincreasing, every $m\leqslant\varrho_t$ then satisfies $g_t(m)\geqslant g_t(\varrho_t)$ and is excluded from $\mathrm{CI}_t$. Thus, using $\inf\varnothing=+\infty$,
  \begin{equation*}
    \mathsf P[\inf\mathrm{CI}_t\geqslant\varrho_t]\geqslant1-\frac{\beta}{2}.
  \end{equation*}

  With $z=\mu-m\geqslant0$, the equation defining $\varrho_t$ reduces to \eqref{root equation I}. The same root bound therefore yields
  \begin{equation*}
    \varrho  _{t} \geqslant \mu-(1+\lambda )\bigg[C_{\alpha}\nu _\alpha(1+u ^{-\alpha })\,\frac{\sum \, \theta _{i}^{1+\alpha }}{\sum \, \theta _{i}}+\frac{L(\beta,\delta)}{\sum \, \theta _{i}}\bigg] .
  \end{equation*}

  Consequently, we derive another one-sided bound, here for $\inf\mathrm{CI}_t$:
  \begin{equation*}
    \mathsf{P}\bigg\{\inf\mathrm{CI}_{t} \geqslant\mu -(1+\lambda )\bigg[C_{\alpha}\nu _\alpha(1+u ^{-\alpha })\,\frac{\sum \, \theta _{i}^{1+\alpha }}{\sum \, \theta _{i}}+\frac{L(\beta,\delta)}{\sum \, \theta _{i}}\bigg]\bigg\} \geqslant 1-\frac{\beta}{2} .
  \end{equation*}

  Combining the two exclusion events gives
  \[
    \mathsf P\{\mathrm{CI}_t\subseteq[\varrho_t,\rho_t]\}\geqslant1-\beta.
  \]

  On this event, the width is at most $\rho_t-\varrho_t$, with width zero if $\mathrm{CI}_t$ is empty. The bounds above consequently imply
  \[
    \mathsf{P}\bigg[|\mathrm{CI}_t|\leqslant\frac{2(1+\lambda)}{\sum \, \theta_i}
    \big((1+u^{-\alpha})C_\alpha\nu_\alpha\sum \, \theta_i^{1+\alpha}+L(\beta,\delta)\big)\bigg]\geqslant1-\beta.
  \]

  This concludes the proof.
\end{proof}

From the proof, we see that the auxiliary parameters $u$ and $\lambda$ arise from two different steps. The parameter $u$ is introduced through the convex splitting in (\ref{convex-splitting-u}), while $\lambda$ quantifies the slack in the upper bound for the smaller positive root of the power equation arising in the width analysis; see (\ref{scaled-root-bound}).
Retaining these proof-generated parameters provides additional flexibility, allowing us to sharpen the width bound in Theorem~\ref{width of CS}. In the special case $\alpha=1$, this leads to the following finite-variance corollary, for which we first introduce a condition (A2*), an analogue of condition (A2).

\textbf{(A2*)} The process $X$ is conditionally square-integrable, in the sense that its conditional second central moment is uniformly bounded by a known constant $\sigma^2>0$:
\begin{equation}\label{second moment assumption}
  \mathsf{E}\big[(X_t-\mu)^2 \mid \mathcal{F}_{t-1}\big] \leqslant \sigma^2  \quad \text{for all } \ t \in \mathbb{N}.
\end{equation}

\begin{corollary}\label{variance CS} 
  Assume that conditions {\rm (A1)} and {\rm (A2*)} hold, and that the coefficients $\{\theta _{t}\}_{t=1}^{\infty }$ are nonrandom and nonnegative. Consider a time $t$ such that $\Sigma_{i=1}^t\,\theta_i>0$. Let $\beta\in(0,1)$, $0<\lambda\leqslant1$, and $0<u<1$ be chosen such that
  \begin{equation*}
    2(1-u  )\frac{\lambda}{(1+\lambda  )^2}\bigg(\frac{\sum_{i=1}^{t} \theta _i}{\sum_{i=1}^{t} \theta _{i}^{2}}\bigg)^2 \geqslant \frac{1+u^{-1}}{2}\, \sigma^2+\frac{L(\beta,\delta)}{\sum_{i=1}^{t} \theta _{i}^2}.
  \end{equation*}

  Then, with probability at least $1-\beta$, the width $|\mathrm{CI}_{t}|$ can be bounded as follows:
  \begin{equation}\label{variance result}
    |\mathrm{CI}_{t}| \leqslant  \frac{2(1+\lambda  )}{\sum_{i=1}^{t} \theta _i} \bigg( \frac{1+u^{-1}}{2}\, \sigma^2 \sum_{i=1}^{t} \theta _{i}^2+L(\beta,\delta)\bigg).
  \end{equation}
\end{corollary}

\begin{remark}[Finite-variance comparison]\label{rem:finite-variance-comparison}
  Similarly, comparing (\ref{variance result}) with the explicit width bound for the Catoni-type CS in \cite[Theorem 10]{CS}, and evaluating the two theorem-level bounds at a common deterministic sequence $\{\theta_t\}_{t=1}^{\infty}$, shows that, for example,  the choice $\lambda=0.5$ and $u=0.75$ yields a sharper explicit upper bound in the finite-variance case as well.
\end{remark}
At this point, one still needs to specify a concrete deterministic sequence $\{\theta_t\}_{t=1}^{\infty}$. In the fixed-time CI setting, \cite{NOC} suggested the choice
\begin{equation*}
  \theta_t = \bigg( \frac{u^\alpha \log ( 2/\delta )}{\alpha \, t \, C_\alpha \,\nu_\alpha} \bigg)^{\frac{1}{1+\alpha}}.
\end{equation*}

Because the explicit width bound derived above for the CS differs from the corresponding fixed-time CI bound, we use in our simulation experiments the alternative tuning
\begin{equation}\label{theta of me}
  \theta_t = \bigg( \frac{\log ( 2/\delta )}{\alpha \,C_\alpha \, \nu_\alpha \,(1 + u^{-\alpha}) t } \bigg)^{\frac{1}{1+\alpha}}.
\end{equation}

The choice in \eqref{theta of me} is a concrete decreasing calibration motivated by the structure of the width bound \eqref{result}. The parameter $\lambda$ does not enter this calibration, but affects both the multiplicative constant in \eqref{result} and the accompanying sufficient condition \eqref{condition}.
For the global tuning sequences listed in Table \ref{table}, the resulting width upper bounds do not attain the iterated-logarithmic shrinkage rate $\mathcal{O}\big(((\log\log t)/t)^{\alpha/(1+\alpha)}\big)$. This comparison motivates the use of the stitched construction, which allows the coefficient choice to vary across different ranges of sample sizes.
\FloatBarrier
\begin{table}[htbp]
  \centering
  \small
  \captionsetup{font=footnotesize, skip=4pt}
  \renewcommand{\arraystretch}{1}
  \setlength{\tabcolsep}{5pt}
  \scalebox{1}{
    \begin{tabular}{|cccc|}
      \hline $\theta_t$                                          & $\Sigma_{i=1}^t\,\theta_i$                                    & $\Sigma_{i=1}^t\,\theta_i^{1+\alpha}$ & Width upper bound                               \\
      \hline $\asymp 1/t$                                        & $\asymp\log t$                                            & $\asymp 1$                        & $1/\log t$                                      \\
      \hline $\asymp ((\log t)/t)^{1/(1+\alpha)}$                & $\asymp (t^\alpha\log t)^{1/(1+\alpha)}$                  & $\asymp\log^2t$                   & $((\log^{2\alpha+1}t)/t^\alpha)^{1/(1+\alpha)}$ \\
      \hline $\asymp (1/t)^{1/(1+\alpha)}$                       & $\asymp t^{\alpha/(1+\alpha)}$                            & $\asymp\log t$                    & $\log t/t^{\alpha/(1+\alpha)}$                  \\
      \hline $\asymp (1/(t\log t))^{1/(1+\alpha)}$               & $\asymp(t^\alpha/\log t)^{1/(1+\alpha)}$                  & $\asymp\log\log t$                & $((\log t)/t^\alpha)^{1/(1+\alpha)}$            \\
      \hline $\asymp [1/(t(\log t)(\log\log t))]^{1/(1+\alpha)}$ & $\asymp [t^\alpha/((\log t)(\log\log t))]^{1/(1+\alpha)}$ & $\asymp\log\log\log t$            & $((\log t)/t^\alpha)^{1/(1+\alpha)}$            \\
      \hline
    \end{tabular}}
  \caption{Rates for the listed global tuning sequences and their cumulative sums, with all other parameters fixed. The final column reports the leading order of the deterministic threshold in \eqref{result}, with multiplicative factors involving $\log\log t$ and further iterated logarithms suppressed only in this column. Powers of $\log t$ and $t$ are retained. The first three columns give the full asymptotic orders. These tuning choices motivate the stitched construction developed below.}
  \label{table}
\end{table}
\FloatBarrier
For $\alpha\in(0,1]$ and $\nu_\alpha>0$, define
\begin{equation}\label{eq:stitched-width-constant}
  A_\alpha=(\mathrm{e}\,\nu_\alpha)^{\frac{1}{1+\alpha}}\,2^{\frac{\alpha}{1+\alpha}}
  (1-\alpha)^{\frac{1-\alpha}{2(1+\alpha)}} \,
  \frac{\mathrm{e}^{\frac{\alpha}{2(1+\alpha)}}(\alpha+\mathrm{e}^{-1})}{\sqrt{1+\alpha}},
\end{equation}
and
\begin{equation}\label{B alpha definition}
  B_{\alpha}(t,\delta) = \bigg(\frac{2\log (\log t+2)+\log (2 / \delta )}{t}\bigg)^{\frac{\alpha }{1+\alpha } }.
\end{equation}

Here, $A_\alpha$ depends only on $\alpha$ and $\nu_\alpha$, whereas
$B_\alpha(t,\delta)$ depends on $\alpha$, $t$, and $\delta$. The definition of $A_\alpha$ at $\alpha=1$ is understood via continuous extension.

\begin{theorem}[Improved stitched $\alpha$-Catoni-type CS]\label{stitching CS}
  Assume {\rm (A1)} and {\rm (A2)}.
  Let $\{\mathrm{CI}_{t}(\theta_*,\delta)\}_{t=1}^{\infty}$ denote the
  $\alpha$-Catoni-type {\rm CS} defined in {\rm (\ref{Catoni-type CS})} with
  constant tuning sequence $\theta_1=\theta_2=\cdots=\theta_*>0$ and error level
  $\delta\in(0,1)$. Let $u\in(0,1)$ and $0<\lambda\leqslant1/\alpha$, and take $C_\alpha$ to be the minimum admissible influence-function constant, with $C_1=1/2$. Define
  \begin{equation}\label{eq:stitched-tuning}
    L_\alpha=(1+u^{-\alpha})C_\alpha\nu_\alpha,\quad t_k=\mathrm{e}^k,\quad
    \delta_k=\frac{\delta}{(k+2)^2},\quad
    \Theta_k=\left(\frac{2\log(2/\delta_k)}{\alpha \, L_\alpha \, t_{k+1}}\right)^{\frac{1}{1+\alpha}},\quad k\in\mathbb N_0.
  \end{equation}

  Then the following stitched sequence,
  \[
    \mathrm{CI}_{t}^{\mathrm{stch}}=\mathrm{CI}_{t}(\Theta_k,\delta_k),
    \qquad t\in[t_k,t_{k+1}),\quad k\in\mathbb{N}_0,
  \]
  is a $(1-\delta)$-{\rm CS} for the unknown mean $\mu$.

  Now, introduce two constants,
  \begin{equation}\label{eq:stitched-constants}
    \kappa=\frac{\pi^2}{3}-\frac52\quad\text{and}\quad
    H=C_\alpha^{-\frac{1}{\alpha}}(1-u)
    \frac{\lambda^{\frac{1}{\alpha}}}{(1+\lambda)^{1+\frac{1}{\alpha}}},
  \end{equation}

  and, for $h>0$, define the following quantities with $k=\lfloor\log t\rfloor$:
  \begin{equation}\label{eq:stitched-threshold}
    \mathcal T(h)=\left\{t\in\mathbb N:t\geqslant3,\ h\,t\geqslant\left(2+\frac{2\,t}{\alpha \,  t_{k+1}}\right)\log\frac2{\delta_k}+\log\kappa\right\}\quad
    \text{and}\quad \tau(h)=\min\mathcal T(h).
  \end{equation}

  We claim that the threshold $\tau(h)$ is finite for every $h>0$ and nonincreasing in $h$. With $B_\alpha(t,\delta)$ as defined in \eqref{B alpha definition}, the stitched sequence satisfies the relation
  \begin{equation}\label{stitching result}
    \mathsf P\left\{\left|\mathrm{CI}_t^{\mathrm{stch}}\right|
    \leqslant2(1+\lambda)(1+u^{-\alpha})^{\frac{1}{1+\alpha}}A_\alpha B_\alpha(t,\delta) \
    \text{ for every integer } \ t\geqslant\tau(H)\right\}\geqslant1-\delta.
  \end{equation}

  The sample-size condition applies only to the width bound, while coverage holds for all $t\in\mathbb N$.
\end{theorem}

\begin{remark}[Width conditions and almost sure shrinkage]\label{rem:width-conditions}
  If $0<u\leqslant u_0$, the time-uniform width bound in Theorem~\ref{stitching CS} holds for every integer
  \[
    t\geqslant\tau\left(\frac{H}{1-u}\right),
  \]
  with the same tuning and displayed width bound. By \eqref{eq:stitched-admissible-times} and the fact that $\tau(h)$ is nonincreasing in $h$, this refinement preserves the width guarantee without increasing the sufficient sample-size threshold $\tau(H)$. In particular, at the benchmark choice $(u,\lambda)=(u_0,\lambda_0)$, the threshold is $\tau(H_{\mathrm{bench}})$, where the subscript ``bench'' denotes this benchmark choice and
  \begin{equation}\label{eq:benchmark-threshold-constant}
    H_{\mathrm{bench}}=\frac{H}{1-u_0}.
  \end{equation}

  To see this, for each deterministic $m\in\mathbb R$ and each $i\geqslant1$, the one-dimensional case of Lemma~8 in Appendix~A.1 of the supplementary material to \cite{WangSGD2021} and (A1)--(A2) yield
  \begin{equation}\label{eq:shifted-moment-bound}
    \mathsf E[|X_i-m|^{1+\alpha}\mid\mathcal F_{i-1}]
    \leqslant4\nu_\alpha+|\mu-m|^{1+\alpha}
    \quad\text{almost surely}.
  \end{equation}

  Since $4=u_0^{-\alpha}\leqslant u^{-\alpha}$, the right-hand side of the last inequality is at most $u^{-\alpha}\nu_\alpha+|\mu-m|^{1+\alpha}$. Thus the proof of Theorem~\ref{stitching CS} applies with the coefficient $(1-u)^{-\alpha}$ of the shifted-moment term replaced by $1$, giving $H/(1-u)$ in place of $H$ while keeping $L_\alpha$, the tuning and the displayed width bound unchanged. At the benchmark choice, $u_0^{-\alpha}=4$, so the central-moment coefficient matches the cited bound exactly.

  For every fixed $\delta\in(0,1)$, with the other parameters fixed as in Theorem~\ref{stitching CS}, the same stitched CS satisfies
  \[
    \left|\mathrm{CI}_t^{\mathrm{stch}}\right|
    =\mathcal{O}\mkern-3mu\left(\left(\frac{\log\log t}{t}\right)^{\frac{\alpha}{1+\alpha}}\right)
    \quad\text{almost surely as } \ t\to\infty.
  \]

  In particular, its width tends to zero almost surely. For completeness, a proof of this almost sure assertion is given in Appendix~\ref{app:almost-sure-shrinkage}.
\end{remark}

The proof of Theorem \ref{stitching CS} relies essentially on an estimate for the upper bound of $|\mathrm{CI}_{t}^{\mathrm{stch}}|$. We need first the following lemma.

\begin{lemma}\label{stitching lemma}
  Suppose that, for $k\in\mathbb{N}_0$, $L_\alpha$, $t_k$, $\delta_k$, and $\Theta_k$ are as defined in \eqref{eq:stitched-tuning}. Let $t$ be any integer and $s$ any real number satisfying
  \[
    t_k\,\mathrm{e}^{\frac{j}{2}}\leqslant s\leqslant t<t_k\, \mathrm{e}^{\frac{j+1}{2}}, \ \mbox{ where } \ j = 0 \mbox{ or } 1.
  \]
  Then we have
  \[
    \frac{L_\alpha s\,\Theta_k^{1+\alpha}+2\log(2/\delta_k)}{s\,\Theta_k}
    \leqslant \mathrm{e}^{\frac{\alpha}{2(1+\alpha)}}(\alpha+\mathrm{e}^{-1})
    \left(\frac2\alpha\right)^{\frac{\alpha}{1+\alpha}}
    (\mathrm{e}\,L_\alpha)^{\frac{1}{1+\alpha}}\,B_\alpha(t,\delta).
  \]
\end{lemma}

\begin{proof}
  Substituting the expression for $\Theta_k$, we have the following chain of relations:
  \begin{align*}
                & \frac{L_\alpha s\,\Theta_k^{1+\alpha}+2\log(2/\delta_k)}{s\,\Theta_k}       \\
    ={}         & \left(\frac2\alpha\right)^{\frac{\alpha}{1+\alpha}}
    L_\alpha^{\frac1{1+\alpha}}
    \left(\frac{\log(2/\delta_k)}s\right)^{\frac{\alpha}{1+\alpha}}
    \left(\frac{s}{t_{k+1}}+\alpha\right)
    \left(\frac{t_{k+1}}s\right)^{\frac1{1+\alpha}}                                           \\
    \leqslant{} & (\alpha+\mathrm{e}^{-1})\left(\frac2\alpha\right)^{\frac{\alpha}{1+\alpha}}
    (\mathrm{e}L_\alpha)^{\frac1{1+\alpha}}
    \left(\frac{\log(2/\delta_k)}s\right)^{\frac{\alpha}{1+\alpha}} \displaybreak[0]          \\
    \leqslant{} & \mathrm{e}^{\frac{\alpha}{2(1+\alpha)}}(\alpha+\mathrm{e}^{-1})
    \left(\frac2\alpha\right)^{\frac{\alpha}{1+\alpha}}
    (\mathrm{e}L_\alpha)^{\frac1{1+\alpha}}
    \left(\frac{\log(2/\delta_k)}t\right)^{\frac{\alpha}{1+\alpha}} \displaybreak[0]          \\
    ={}         & \mathrm{e}^{\frac{\alpha}{2(1+\alpha)}}(\alpha+\mathrm{e}^{-1})
    \left(\frac2\alpha\right)^{\frac{\alpha}{1+\alpha}}
    (\mathrm{e}L_\alpha)^{\frac1{1+\alpha}}
    \left(\frac{\log(2/\delta)+2\log(k+2)}t\right)^{\frac{\alpha}{1+\alpha}} \displaybreak[0] \\
    \leqslant{} & \mathrm{e}^{\frac{\alpha}{2(1+\alpha)}}(\alpha+\mathrm{e}^{-1})
    \left(\frac2\alpha\right)^{\frac{\alpha}{1+\alpha}}
    (\mathrm{e}L_\alpha)^{\frac1{1+\alpha}}
    \left(\frac{\log(2/\delta)+2\log(\log t+2)}t\right)^{\frac{\alpha}{1+\alpha}}             \\
    ={}         & \mathrm{e}^{\frac{\alpha}{2(1+\alpha)}}(\alpha+\mathrm{e}^{-1})
    \left(\frac2\alpha\right)^{\frac{\alpha}{1+\alpha}}
    (\mathrm{e}L_\alpha)^{\frac1{1+\alpha}}\,B_\alpha(t,\delta).
  \end{align*}

  For the first inequality, note that $\mathrm{e}^{-1}\leqslant s/t_{k+1}<1$ and
  \[
    \frac{d}{dx}\left[(x+\alpha)x^{-\frac1{1+\alpha}}\right]
    =\frac{\alpha(x-1)}{1+\alpha}\,
    x^{-\frac{2+\alpha}{1+\alpha}}\leqslant0,\qquad 0<x\leqslant1.
  \]

  Hence $(s/t_{k+1}+\alpha)(t_{k+1}/s)^{1/(1+\alpha)}\leqslant(\alpha+\mathrm{e}^{-1})\mathrm{e}^{1/(1+\alpha)}$.
  This proves the lemma. The argument uses only $L_\alpha>0$ and therefore applies to any positive value of $L_\alpha$, with $\Theta_k$ defined accordingly.
\end{proof}

\begin{proof}[Proof of Theorem~\ref{stitching CS}]
  We first establish coverage. For each $k\in\mathbb N_0$, Theorem~\ref{Catoni-type confidence sequence} shows that $\{\mathrm{CI}_t(\Theta_k,\delta_k)\}_{t=1}^{\infty}$ is a $(1-\delta_k)$-CS for $\mu$. The time blocks $T_k=\mathbb N\cap[t_k,t_{k+1})$, $k\in\mathbb N_0$, form a partition of $\mathbb N$, and $\sum_{k=0}^{\infty}\,\delta_k=\delta\sum_{k=0}^{\infty}\,\frac1{(k+2)^2}<\delta$.

  Applying Lemma~\ref{Stitching method} with this partition and the sequence $\{\delta_k\}_{k=0}^{\infty}$, we conclude that the stitched sequence
  \[
    \mathrm{CI}_t^{\mathrm{stch}}=\mathrm{CI}_t(\Theta_k,\delta_k),
    \qquad t\in T_k,\quad k\in\mathbb N_0,
  \]
  is a $(1-\delta)$-CS for $\mu$, giving the claimed coverage for all $t\in\mathbb N$.

  For the width bound, fix $k\in\mathbb N$ and $j\in\{0,1\}$, and let $I_{k,j}:=[t_k\mathrm{e}^{j/2},t_k\mathrm{e}^{(j+1)/2})$. Using \eqref{eq:log-error-budget} and \eqref{eq:stitched-constants}, set
  \begin{equation}\label{eq:stitched-proof-errors}
    \varepsilon_k=\frac{\delta_k}{\kappa},\qquad
    L(\varepsilon_k,\delta_k)=2\log\frac2{\delta_k}+\log\kappa.
  \end{equation}

  We use $g_t,B_t^\pm,K_t,J_t$ defined in \eqref{eq:width-score}, \eqref{eq:width-envelopes} and \eqref{eq:width-root-coefficients}, with $\theta_i=\Theta_k$ for every $i$, $\beta=\varepsilon_k$ and $\delta$ replaced by $\delta_k$. If $I_{k,j}\cap\mathcal T(H)=\varnothing$, we set $E_{k,j}=\varnothing$. Otherwise, the threshold condition gives $H>L_\alpha\Theta_k^{1+\alpha}$, and the deterministic real number
  \[
    s=\max\left\{t_k\mathrm{e}^{{j}/{2}}, \
    \frac{L(\varepsilon_k,\delta_k)}
    {H-L_\alpha\Theta_k^{1+\alpha}}\right\}
  \]
  satisfies $t_k\mathrm{e}^{j/2}\leqslant s\leqslant t<t_k\mathrm{e}^{(j+1)/2}$ for every $t\in I_{k,j}\cap\mathcal T(H)$. Evaluating the constant-weight coefficient formulas at this real $s$, we have
  \[
    K_s=C_\alpha(1-u)^{-\alpha}\Theta_k^\alpha,\qquad
    J_s=L_\alpha\Theta_k^\alpha+
    \frac{L(\varepsilon_k,\delta_k)}{s\,\Theta_k},
  \]
  and $\Theta_kJ_s\leqslant H$. Thus the deterministic comparison calculation following \eqref{parameter condition} in the proof of Theorem~\ref{width of CS} applies. Take the fixed comparison points
  \[
    \rho=\mu+(1+\lambda)J_s,\qquad
    \varrho=\mu-(1+\lambda)J_s.
  \]

  Since $K_t=K_s$ and $J_t\leqslant J_s$ for $t\geqslant s$, that calculation gives, for every integer $t\geqslant s$,
  \[
    B_t^+(\rho)=-B_t^-(\varrho)
    \leqslant-C_\alpha\nu_\alpha \,t\,\Theta_k^{1+\alpha}-\log\frac2{\delta_k}.
  \]

  Define the bound-crossing event
  \begin{equation}\label{eq:stitched-failure-events}
    E_{k,j}
    =\bigl\{\exists\,t\in\mathbb N:g_t(\rho)\geqslant B_t^+(\rho)\bigr\}
    \cup\bigl\{\exists\,t\in\mathbb N:g_t(\varrho)\leqslant B_t^-(\varrho)\bigr\}.
  \end{equation}

  Applying Ville's inequality to the two supermartingales $\widetilde M_t^+(\rho,u)$ and $\widetilde M_t^-(\varrho,u)$ defined in \eqref{eq:compensated-supermartingales}, with boundary $2/\varepsilon_k$ for each, gives $\mathsf P(E_{k,j})\leqslant\varepsilon_k$. Both processes use the constant weight $\Theta_k$ and all observations starting from $X_1$. On $E_{k,j}^c$, the strict bounds $g_t(\rho)<B_t^+(\rho)$ and $g_t(\varrho)>B_t^-(\varrho)$ hold for all $t\in\mathbb N$. Applying the monotone-exclusion argument from the proof of Theorem~\ref{width of CS}, we obtain
  \[
    \mathrm{CI}_t^{\mathrm{stch}}\subset(\varrho,\rho)
    \quad\text{for all } \ t\in I_{k,j}\cap\mathcal T(H).
  \]

  Using $L(\varepsilon_k,\delta_k)\leqslant2\log(2/\delta_k)$ and Lemma~\ref{stitching lemma}, on this same event we have
  \[
    \begin{aligned}
      \left|\mathrm{CI}_t^{\mathrm{stch}}\right| \leqslant{} & 2(1+\lambda)J_s                                \\
      \leqslant{}                                            & 2(1+\lambda)
      \frac{L_\alpha s\,\Theta_k^{1+\alpha}+2\log(2/\delta_k)}{s\,\Theta_k}                                   \\
      \leqslant{}                                            & 2(1+\lambda)(1+u^{-\alpha})^{\frac1{1+\alpha}}
      A_\alpha B_\alpha(t,\delta)
    \end{aligned}
  \]
  simultaneously for all such $t$.

  Taking a union bound over the two time subintervals for every $k\geqslant1$ gives
  \begin{equation}\label{eq:stitched-error-budget}
    \mathsf P\left(\bigcup_{k\geqslant1}\bigcup_{j=0}^1 E_{k,j}\right)
    \leqslant\sum_{k\geqslant1}\sum_{j=0}^1\mathsf P(E_{k,j})
    \leqslant\frac2\kappa\sum_{k\geqslant1}\delta_k
    =\frac{2\delta}{\kappa}\sum_{n\geqslant3}\frac1{n^2}
    =\delta.
  \end{equation}

  Hence the width bound holds simultaneously for all $t\in\mathcal T(H)$ with probability at least $1-\delta$.

  It remains to verify the properties of the threshold. Fix $h>0$. For an integer $t\geqslant3$, with $k=\lfloor\log t\rfloor$, the condition in \eqref{eq:stitched-threshold} is equivalent to
  \[
    h\geqslant\frac{2\log(2/\delta_k)}{\alpha t_{k+1}}
    +\frac{2\log(2/\delta_k)+\log\kappa}{t}.
  \]

  Since $\log(2/\delta_k)=\log(2/\delta)+2\log(k+2)$ and $t\geqslant \mathrm{e}^k$, the right-hand side tends to zero as $t\to\infty$. Hence $\mathcal T(h)$ contains all sufficiently large integers, so it is nonempty and $\tau(h)<\infty$.

  We next show that the condition continues to hold once it is satisfied. For this purpose, extend the deterministic expression on the right to real $t\geqslant \mathrm{e}$, with $k=\lfloor\log t\rfloor$. Within the fixed time interval $[t_k,t_{k+1})$, both $\delta_k$ and $t_{k+1}$ are constant. The first term is therefore constant, and the second decreases with $t$ because $2\log(2/\delta_k)+\log\kappa>0$.

  At the boundary $t_{k+1}$, the left-hand limit from $[t_k,t_{k+1})$ exceeds the value at the beginning of the next interval $[t_{k+1},t_{k+2})$:
  \[
    \frac{\left(2+\frac2\alpha\right)\log(2/\delta_k)+\log\kappa}{t_{k+1}}>\frac{\left(2+\frac2{\alpha \mathrm{e}}\right)\log(2/\delta_{k+1})+\log\kappa}{t_{k+1}}.
  \]

  Indeed, $t_{k+2}=\mathrm{e}\,t_{k+1}$, and
  \[
    \frac{\log(2/\delta_{k+1})}{\log(2/\delta_k)}
    \leqslant\frac{\log18}{\log8}
    <\frac{\alpha+1}{\alpha+\mathrm{e}^{-1}}.
  \]

  Thus the right-hand side decreases both within these time intervals and across their boundaries. In particular, every integer after the first qualifying time also satisfies the condition, giving
  \begin{equation}\label{eq:stitched-admissible-times}
    \mathcal T(h)=\{t\in\mathbb N:t\geqslant\tau(h)\}.
  \end{equation}

  Finally, for $0<h_1\leqslant h_2$, the definition implies $\mathcal T(h_1)\subseteq\mathcal T(h_2)$, and hence $\tau(h_2)\leqslant\tau(h_1)$.

  Taking $h=H$, the width bound already established on $\mathcal T(H)$ therefore holds simultaneously for every integer $t\geqslant\tau(H)$ with probability at least $1-\delta$. This completes the proof.
\end{proof}

Table~\ref{table:sufficient-sample-sizes} reports the permanent thresholds $\tau(H_{\mathrm{bench}})$ for the time-uniform width bound in Remark~\ref{rem:width-conditions} under the benchmark choice $(u,\lambda)=(u_0,\lambda_0)$. Here, ``permanent'' means that the sufficient condition holds for every integer $t$ at or above the threshold, as shown in \eqref{eq:stitched-admissible-times}. These thresholds are derived from sufficient conditions in the proof and may therefore be conservative. Although $\tau(h)$ is the smallest integer satisfying the condition in \eqref{eq:stitched-threshold}, we do not claim that it is the smallest possible threshold for the time-uniform width guarantee.

\FloatBarrier
\begin{table}[htbp]
  \centering
  \small
  \captionsetup{font=footnotesize, skip=4pt}
  \renewcommand{\arraystretch}{1}
  \scalebox{1}{
    \begin{tabular}{c|cccc}
      \hline
      $\alpha \backslash \delta$ & $0.1$ & $0.05$ & $0.01$ & $0.001$ \\
      \hline
      $1$                        & $45$  & $55$   & $55$   & $64$    \\
      $1/2$                      & $63$  & $76$   & $121$  & $149$   \\
      $1/3$                      & $404$ & $404$  & $445$  & $792$   \\
      \hline
    \end{tabular}}
  \caption{Permanent thresholds $\tau(H_{\mathrm{bench}})$ for time-uniform width control under the benchmark choice $(u,\lambda)=(u_0,\lambda_0)$, using the refined condition in Remark~\ref{rem:width-conditions}. These sufficient thresholds do not identify empirical onset times, and coverage holds from $t=1$.}
  \label{table:sufficient-sample-sizes}
\end{table}
\FloatBarrier

\begin{remark}[Joint dependence on time and confidence level]\label{joint dependence}
  The time-uniform width bound in Theorem~\ref{stitching CS} attains the iterated-logarithmic shrinkage rate for fixed $\delta$, improving the asymptotic order of the width bounds associated with the global tuning sequences listed in Table~\ref{table}. Motivated by the joint asymptotic perspective as described in \cite{DBJ2026} for fixed-time CIs, we fix $\alpha$, $\nu_\alpha$, $u$, and $\lambda$ as in Theorem~\ref{stitching CS}, and see that the dependence on $t$ and $\delta$ is entirely through $B_\alpha(t,\delta)$.

  For fixed $\delta$, \eqref{B alpha definition} gives
  \[
    B_\alpha(t,\delta)\asymp
    \left(\frac{\log\log t}{t}\right)^{\frac{\alpha}{1+\alpha}},
    \quad \mbox { as } \ t\to\infty.
  \]

  For fixed $t$, $B_\alpha(t,\delta)$ increases as $\delta$ decreases and, as a
  deterministic expression, satisfies
  \[
    B_\alpha(t,\delta)\asymp
    \left(\frac{\log(1/\delta)}{t}\right)^{\frac{\alpha}{1+\alpha}},
    \quad \mbox{ as } \ \delta\to0.
  \]

  The condition $t\geqslant\tau(H)$, as characterized by \eqref{eq:stitched-admissible-times} with $h=H$, specifies a range of $\delta$ over which the high-probability time-uniform width bound in \eqref{stitching result} is guaranteed.

  To study the joint regime, let $\{t_n\}_{n=1}^{\infty}\subset\mathbb N$ be an increasing sequence with $t_n\to\infty$, and let $\{\delta_n\}_{n=1}^{\infty}\subset(0,1)$ be a decreasing sequence with $\delta_n\to0$. Here and below, any asymptotic relation is understood to hold as $n\to\infty$, the latter, `$n \to \infty$', being omitted just for simplicity of writing. Set
  \begin{equation}\label{eq:joint-regime-notation}
    T_n=2\log(\log t_n+2),\qquad
    D_n=\log(2/\delta_n),\qquad k_n=\lfloor \log t_n \rfloor.
  \end{equation}

  Then
  \[
    B_\alpha(t_n,\delta_n)
    =\left(\frac{T_n+D_n}{t_n}\right)^{\frac{\alpha}{1+\alpha}}.
  \]

  Since $T_n=o(t_n)$, we have
  \[
    B_\alpha(t_n,\delta_n)\rightarrow 0
    \quad\Longleftrightarrow\quad D_n=o(t_n).
  \]

  If $D_n=\mathcal{O}(T_n)$, then
  \[
    B_\alpha(t_n,\delta_n)
    =\widetilde{\mathcal{O}}\left(
    \left(\frac{\log\log t_n}{t_n}\right)^{\frac{\alpha}{1+\alpha}}
    \right).
  \]

  Thus the iterated-logarithmic order is preserved, including when $\delta_n=(\log
    t_n)^{-r}$ for fixed $r>0$. The first transition scale is $\log(1/\delta_n)=\widetilde{\mathcal{O}}(\log\log t_n)$, at which the time-uniform term $T_n$ and the confidence-level term $D_n$ are comparable. Equivalently, for all sufficiently large $n$,
  \[
    \delta_n=(\log t_n)^{-r_n}
    \quad\Longleftrightarrow\quad
    t_n=\exp\left(\delta_n^{-1/r_n}\right),
  \]
  where $r_n:=\log(1/\delta_n)/\log\log t_n$ is bounded away from zero and
  infinity.

  If instead, $T_n=o(D_n)$ and $D_n=o(t_n)$, the confidence-level term $D_n$ dominates and
  \[
    B_\alpha(t_n,\delta_n)
    \asymp\left(\frac{\log(1/\delta_n)}{t_n}\right)^{\frac{\alpha}{1+\alpha}}
    \longrightarrow \ 0.
  \]

  For example, $\delta_n=t_n^{-r}$ with fixed real $r>0$ gives
  \[
    B_\alpha(t_n,\delta_n)
    \asymp\left(\frac{\log t_n}{t_n}\right)^{\frac{\alpha}{1+\alpha}}
    \longrightarrow \ 0.
  \]

  In both regimes, see above, and more generally whenever $D_n=o(t_n)$, the inequality in \eqref{eq:stitched-threshold}, with $h=H$ and $(t,\delta)=(t_n,\delta_n)$, eventually holds, since its right-hand side is $\widetilde{\mathcal O}(T_n+D_n)$. Hence \eqref{stitching result} and $\delta_n\to0$ imply that
  \[
    \bigl|\mathrm{CI}_{t_n}^{\mathrm{stch}}(\delta_n)\bigr| \
    \xrightarrow{\mathsf P} \ 0,
  \]
  where $\mathrm{CI}_t^{\mathrm{stch}}(\delta_n)$ denotes the stitched CS constructed at error level $\delta_n$. This is a sufficient condition for the random widths to vanish; the equivalence $B_\alpha(t_n,\delta_n)\to0\Longleftrightarrow D_n=o(t_n)$ concerns the deterministic width bound.

  If, in addition, $\Sigma_{n=1}^{\infty}\,\delta_n<\infty$, the widths converge to zero almost surely along these deterministic sequences on the same data stream:
  \[
    \bigl|\mathrm{CI}_{t_n}^{\mathrm{stch}}(\delta_n)\bigr| \
    \xrightarrow{\mathrm{a.s.}} \ 0.
  \]

  Here the error level varies with $n$, whereas Remark~\ref{rem:width-conditions} concerns a single CS at a fixed error level. The proof of this is given in Appendix~\ref{app:almost-sure-shrinkage}.

  Finally, if $D_n/t_n \ \to \ \gamma$ for some real $\gamma \in(0,\infty)$, then
  \[
    B_\alpha(t_n,\delta_n) \ \longrightarrow \
    \gamma^{\frac{\alpha}{1+\alpha}}.
  \] \

  Whenever the width condition holds for all sufficiently large $n$, the bound in \eqref{stitching result} is of constant order. If $D_n/t_n\to\infty$, this condition eventually fails, and Theorem~\ref{stitching CS} no longer guarantees the stated high-probability time-uniform width bound. The two transition scales, $\log(1/\delta_n)=\widetilde{\mathcal{O}}(\log\log t_n)$ and $\log(1/\delta_n)=\widetilde{\mathcal{O}}(t_n)$, therefore describe the balance between the two terms and the loss of a vanishing bound, respectively.
\end{remark}

Having described the asymptotic dependence on $t$ and $\delta$, we next
compare the multiplicative constants in the stitched time-uniform width bounds.
For this purpose, Theorem~\ref{stitching CS of Ramdas} below gives the stitched
version of the $\alpha$-Catoni-type CS in \cite[Theorem~12]{CS}, obtained using
the same stitching scheme. The superscripts $(\mathrm{WR})$ and
$(\mathrm{C})$ refer to Wang--Ramdas and Chen et al., respectively.
The time-uniform width bound is verified directly under the sufficient condition below, using the same simultaneous-width arguments as in Theorem~\ref{stitching CS}; see Appendix~\ref{app:wr-stitched}.

\begin{theorem}[Stitched $\alpha$-Catoni-type CS]\label{stitching CS of Ramdas}
  Assume {\rm (A1)} and {\rm (A2)}.
  Let $\psi_\alpha^{\mathrm{(C)}}$ be a nondecreasing influence function satisfying \eqref{influence function of Chen}. Define $\{\mathrm{CI}_t^{\rm(WR)}(\theta_*,\delta)\}_{t=1}^{\infty}$ as the following $\alpha$-Catoni-type {\rm CS}:
  \begin{equation}\label{Catoni-type CS of Ramdas}
    \mathrm{CI}_t^{\rm(WR)}(\{\theta_t\},\delta)
    =\left\{m\in\mathbb R:
    \left|\sum_{i=1}^t\psi_\alpha^{\mathrm{(C)}}\bigl(\theta_i(X_i-m)\bigr)\right|
    \leqslant\frac{\nu_\alpha}{1+\alpha}\sum_{i=1}^t\theta_i^{1+\alpha}
    +\log\frac2\delta\right\},
  \end{equation}
  with a constant tuning sequence $\theta_1=\theta_2=\cdots=\theta_*>0$ and error level $\delta\in(0,1)$. Define
  \begin{equation}\label{eq:wr-stitched-tuning}
    L_\alpha^{\rm(WR)}=\frac{5\nu_\alpha}{1+\alpha},\quad t_k=\mathrm{e}^k,\quad
    \delta_k=\frac{\delta}{(k+2)^2},\quad
    \Theta_k^{\rm(WR)}=\left(\frac{2\log(2/\delta_k)}
    {\alpha L_\alpha^{\rm(WR)}t_{k+1}}\right)^{\frac{1}{1+\alpha}},\quad
    k\in\mathbb N_0.
  \end{equation}

  Then the following stitched sequence,
  \[
    \mathrm{CI}_t^{\mathrm{stch(WR)}}
    =\mathrm{CI}_t^{\rm(WR)}(\Theta_k^{\rm(WR)},\delta_k),
    \quad t\in[t_k,t_{k+1}),\quad k\in\mathbb N_0,
  \]
  is a $(1-\delta)$-{\rm CS} for the unknown mean $\mu$.

  Further, put
  \begin{equation}\label{eq:wr-threshold-constant}
    H_{\mathrm{WR}}=\frac{\alpha^{1/\alpha}}{1+\alpha},
  \end{equation}
  and let $\tau$ be defined in \eqref{eq:stitched-threshold}. With $B_\alpha(t,\delta)$ as defined in \eqref{B alpha definition}, the stitched sequence satisfies the following bound:
  \begin{equation}\label{stitching result of Ramdas}
    \mathsf P\left\{\left|\mathrm{CI}_t^{\mathrm{stch(WR)}}\right|
    \leqslant2A_\alpha^{\rm(WR)}B_\alpha(t,\delta) \
    \text{ for every integer } \ t\geqslant\tau(H_{\mathrm{WR}})\right\}\geqslant1-\delta.
  \end{equation}
  Here
  \begin{equation}\label{eq:wr-width-constant}
    A_\alpha^{\rm(WR)}=(5\,\mathrm{e}\,\nu_\alpha)^{\frac{1}{1+\alpha}}
    \left(\frac2\alpha\right)^{\frac{\alpha}{1+\alpha}}
    (1+\alpha)^{\frac{\alpha}{1+\alpha}}\,
    \mathrm{e}^{\frac{\alpha}{2(1+\alpha)}}\,(\alpha+\mathrm{e}^{-1}),
  \end{equation}
  a constant depending only on $\alpha$ and $\nu_\alpha$.

  The sample-size condition applies only to the width bound, while coverage holds for all $t\in\mathbb N$.
\end{theorem}

\begin{remark}[Comparison of stitched time-uniform width bounds]\label{choose of omega in SM}
  The bounds in \eqref{stitching result} and
  \eqref{stitching result of Ramdas} share the factor $B_\alpha(t,\delta)$
  and hence have the same asymptotic dependence in $t$ and $\delta$.
  Under the benchmark choice $(u,\lambda)=(u_0,\lambda_0)$, a direct comparison using \eqref{eq:stitched-width-constant} and \eqref{eq:wr-width-constant} shows that the improved stitched $\alpha$-Catoni-type CS in Theorem~\ref{stitching CS} has a strictly smaller constant in its width bound than the stitched CS in Theorem~\ref{stitching CS of Ramdas} for $0<\alpha<1$.

  More explicitly, at these benchmark parameters and for $0<\alpha<1$, the ratios of the displayed width constants and the sufficient-condition constants in \eqref{eq:benchmark-threshold-constant} and \eqref{eq:wr-threshold-constant} satisfy the relations
  \[
    \frac{(1+\lambda_0)(1+u_0^{-\alpha})^{\frac{1}{1+\alpha}}A_\alpha}{A_\alpha^{\rm(WR)}}
    =((1+\alpha)C_\alpha)^{\frac{1}{1+\alpha}}<1 \quad \mbox{and} \quad
    \frac{H_{\mathrm{bench}}}{H_{\mathrm{WR}}}
    =((1+\alpha)C_\alpha)^{-\frac{1}{\alpha}}>1.
  \]

  Thus, for $0<\alpha<1$, the sufficient condition for the proposed benchmark bound holds whenever that for the Wang--Ramdas bound holds. At $\alpha=1$, these benchmark bounds and the sufficient conditions coincide.
\end{remark}

%% file: section/Simulationresults.tex
\section{Simulation study}\label{experiments}
This section examines both the efficiency and the empirical validity of the proposed CSs. In particular, we compare interval widths over time and empirical miscoverage for two heavy-tailed models: a centered Pareto distribution with shape parameter $1.8$, whose density is
\[
  f_{\mathrm P}(x)=
  \begin{cases}
    1.8\left(x+\dfrac94\right)^{-2.8}, & x>-\dfrac54,          \\[8pt]
    0,                                 & x\leqslant -\dfrac54,
  \end{cases}
\]
and a centered Student's $\rm t$-distribution with $2$ degrees of freedom, whose density is
\[
  f_{\mathrm t}(x)=\Bigl(x^2+2\Bigr)^{-\frac{3}{2}},\quad x\in\mathbb{R}.
\]

For both models, the second-order moment, and hence the variance, is infinite.
However, a random variable with density $f_{\mathrm P}$ has a finite absolute central moment of order $1+\alpha$ for $0<\alpha<0.8$, whereas one with density $f_{\mathrm t}$ has such a moment for $0<\alpha<1$.
Accordingly, throughout the simulation study we set $\alpha=0.5$. The a priori bound $\nu_\alpha$ is taken to be $4$ in all experiments, based on a calculation from the above densities.
All observations are sampled i.i.d., so these experiments illustrate this subclass of (A1)--(A2).

For both models, in each of the two procedures, non-stitched and stitched, we use the influence function
\begin{equation}\label{simulation-influence-function}
  \psi_\alpha(x)=
  \begin{cases}
    \log\bigl(1+x+C_\alpha x^{1+\alpha}\bigr),   & x>0,         \\[3pt]
    -\log\bigl(1-x+C_\alpha|x|^{1+\alpha}\bigr), & x\leqslant0.
  \end{cases}
\end{equation}

The Wang--Ramdas comparison uses the same form as \eqref{simulation-influence-function} with $C_\alpha$ replaced by $C=1/(1+\alpha)$. Both influence functions satisfy their respective bounds and are continuous, strictly increasing, and unbounded in both directions. Hence, with positive tuning parameters, the estimating and interval-endpoint equations have unique solutions and can be solved numerically.

Across methods, all procedures are supplied with the same values of $\alpha$, $\nu_\alpha$, and $\delta$; no distribution-specific retuning is performed. For the non-stitched procedures, we use the quantities
\begin{equation}\label{eq:simulation-nonstitched-tuning}
  \theta_t =\bigg( \frac{\log ( 2/\delta )}{\alpha \,C_\alpha\, \nu_\alpha\,(1 + u^{-\alpha}) \, t } \bigg)^{\frac{1}{1+\alpha}}
  \qquad\text{and}\qquad
  \theta_t^{\rm (WR)} =\frac{1}{2}\bigg( \frac{ 2(1+\alpha) \log ( 2/\delta )}{ t \, \nu_\alpha } \bigg)^{\frac{1}{1+\alpha}},
\end{equation}
where the proposed procedure is based on the benchmark value $u=u_0$, as in Remark~\ref{choose of omega}. The sequence $\{\theta_t\}$ is chosen according to \eqref{theta of me}, whereas $\{\theta_t^{\rm (WR)}\}$ is the deterministic calibration used in \cite[Theorem 12]{CS}. We then construct the corresponding confidence sequences $\{\mathrm{CI}_t\}_{t=1}^{\infty}$ from \eqref{Catoni-type CS} and $\{\mathrm{CI}_t^{\mathrm{(WR)}}\}_{t=1}^{\infty}$ from \eqref{Catoni-type CS of Ramdas}. The comparison concerns the complete calibrated procedures, including their influence functions and respective deterministic tuning sequences.

\begin{figure}[htbp]
  \centering
  \captionsetup{size=footnotesize, skip=4pt}
  \begin{subfigure}[b]{0.495\textwidth}
    \includegraphics[width=\textwidth]{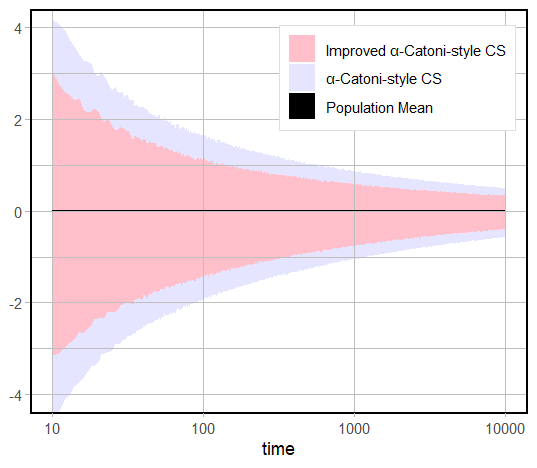}
    \caption*{Centered Pareto$(1.8)$ model}
  \end{subfigure}
  \hfill
  \begin{subfigure}[b]{0.495\textwidth}
    \includegraphics[width=\textwidth]{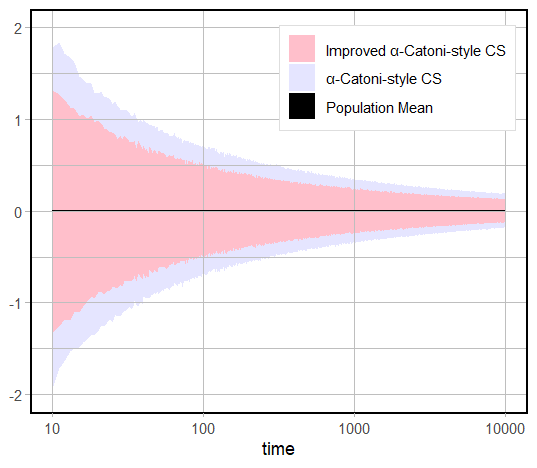}
    \caption*{Centered $t_2$ model}
  \end{subfigure}
  \caption{Non-stitched procedures at level $\delta=0.05$ with $\alpha=0.5$. The lower and upper boundaries of each band are the empirical means of the corresponding interval endpoints over $200$ independent runs. Their vertical separation is the empirical mean interval width. In both heavy-tailed models, the proposed $\alpha$-Catoni-type CS has a smaller empirical mean interval width than the Wang--Ramdas $\alpha$-Catoni-type CS over the displayed horizon.}
  \label{fig:nonstitched-width}
\end{figure}

We first compare the non-stitched procedures at error level $\delta=0.05$, with emphasis on the shrinkage rate behavior  as $t\to\infty$. Figure \ref{fig:nonstitched-width} plots, for each $t$, the empirical mean lower and upper interval endpoints based on $200$ independent Monte Carlo runs for the improved procedure and for the Wang--Ramdas procedure. The vertical separation of these mean endpoints equals the empirical mean interval width.
In both models, the proposed $\alpha$-Catoni-type CS has a smaller empirical mean interval width than the Wang--Ramdas CS over the displayed time range.

To complement the width comparison, we next examine empirical cumulative miscoverage for the two non-stitched procedures considered above. For a given monitoring horizon $t$, this quantity is the proportion of runs in which the interval sequence $\{\mathrm{CI}_s\}_{1 \leqslant s \leqslant t}$ fails to cover $\mu$ at least once up to time $t$.
The estimates are based on $4000$ independent repetitions, so that the corresponding Monte Carlo error is small relative to the nominal level $\delta=0.05$.

\FloatBarrier

\FloatBarrier

Under the centered Pareto$(1.8)$ model, the cumulative miscoverage curves are monotone and nearly degenerate, so we summarize them in Table \ref{tab:pareto-cum-miscoverage} through the first miscoverage times. For the proposed procedure, only $5$ of the $4000$ runs exhibit a failure to cover $\mu$ over the displayed monitoring horizon, and the first failure in each of these runs occurs at one of the first two monitoring times.
For the Wang--Ramdas procedure, no failure is observed in these runs. Thus, in this example, the cumulative miscoverage of both procedures remains far below the nominal level over the displayed horizon.

\FloatBarrier
\begin{table}[htbp]
  \centering
  \small
  \captionsetup{font=footnotesize, skip=4pt}
  \renewcommand{\arraystretch}{1}
  \scalebox{1}{
    \begin{tabular}{c|c|c|c}
      \hline
      First miscoverage time           & $t=1$ & $t=2$ & No miscoverage \\
      \hline
      Improved $\alpha$-Catoni-type CS & 4     & 1     & 3995           \\
      $\alpha$-Catoni-type CS of WR    & 0     & 0     & 4000           \\
      \hline
    \end{tabular}
  }
  \caption{Centered Pareto$(1.8)$ model with $\alpha=0.5$ and $\delta=0.05$. The table reports the empirical distribution of the first miscoverage time over $4000$ independent runs. For both non-stitched procedures, the observed cumulative miscoverage is concentrated at the earliest monitoring times and remains far below the nominal level over the displayed horizon.}
  \label{tab:pareto-cum-miscoverage}
\end{table}
\FloatBarrier

\begin{figure}[htbp]
  \centering
  \captionsetup{size=footnotesize, skip=4pt}
  \includegraphics[width=0.7\textwidth]{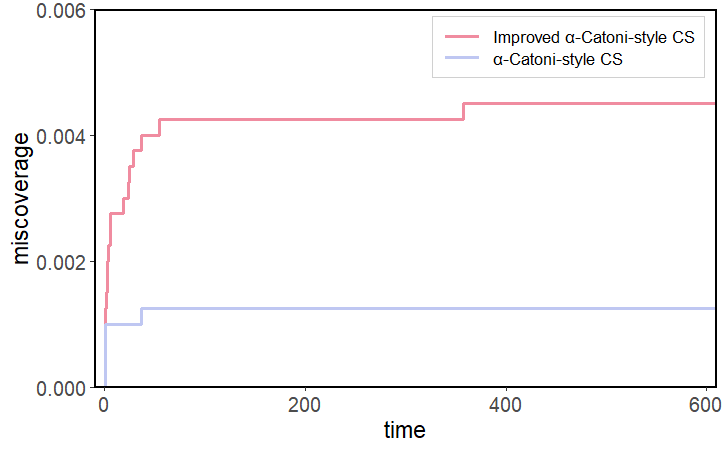}
  \caption{Empirical cumulative miscoverage rates under the centered $t_2$ model for the proposed $\alpha$-Catoni-type CS and the Wang--Ramdas $\alpha$-Catoni-type CS, with $\alpha=0.5$ and $\delta=0.05$. The curves are based on $4000$ independent runs. Both procedures remain well below the nominal miscoverage level over the displayed monitoring horizon.}
  \label{fig:t-cum-miscoverage}
\end{figure}

Under the centered $t_2$ model, the cumulative miscoverage curves are more informative and are therefore shown in Figure \ref{fig:t-cum-miscoverage}. The cumulative miscoverage rates of both procedures increase stepwise under continuous monitoring, but they remain far below the nominal level throughout the displayed horizon. Thus, in this second example as well, the smaller empirical mean interval widths produced by the proposed method are not obtained by materially degrading empirical validity.

\FloatBarrier

Next we examine how the empirical mean interval widths vary with the confidence level at fixed sample sizes. Fixing $t\in\{100,1000,10000\}$ and varying $\delta$, we use the same deterministic sequences $\{\theta_t\}$ and $\{\theta_t^{\rm (WR)}\}$ from \eqref{eq:simulation-nonstitched-tuning} and estimate the mean interval widths from $200$ independent runs for each configuration. The resulting comparisons under the two heavy-tailed models are reported together in Table \ref{delta}.

\FloatBarrier
\begin{table}[htbp]
  \centering
  \small
  \captionsetup{font=footnotesize, skip=4pt}
  \renewcommand{\arraystretch}{1}
  \scalebox{1}{
    \begin{tabular}{lcccccc}
      \toprule
               & \multicolumn{3}{c}{Improved $\alpha$-Catoni-type CS} & \multicolumn{3}{c}{$\alpha$-Catoni-type CS of WR}                                              \\
      \cmidrule(lr){2-4} \cmidrule(lr){5-7}
      $\delta$ & $t=100$                                              & $t=1000$                                          & $t=10000$ & $t=100$ & $t=1000$ & $t=10000$ \\
      \midrule
      0.2      & 2.12919                                              & 1.13932                                           & 0.62651   & 2.94888 & 1.62607  & 0.91268   \\
      0.1      & 2.34596                                              & 1.24633                                           & 0.68367   & 3.26571 & 1.78010  & 0.99643   \\
      0.05     & 2.54187                                              & 1.33897                                           & 0.73279   & 3.55163 & 1.91267  & 1.06812   \\
      0.02     & 2.76818                                              & 1.44595                                           & 0.78922   & 3.89399 & 2.06854  & 1.15076   \\
      0.01     & 2.93367                                              & 1.51862                                           & 0.82752   & 4.14518 & 2.17451  & 1.20599   \\
      0.005    & 3.08924                                              & 1.58578                                           & 0.86202   & 4.38471 & 2.26966  & 1.25758   \\
      0.002    & 3.27635                                              & 1.66688                                           & 0.90462   & 4.68365 & 2.39164  & 1.31923   \\
      0.001    & 3.41521                                              & 1.72451                                           & 0.93431   & 4.91226 & 2.47526  & 1.36287   \\
      \bottomrule
    \end{tabular}
  }

  \vspace{0.8em}

  \scalebox{1}{
    \begin{tabular}{lcccccc}
      \toprule
               & \multicolumn{3}{c}{Improved $\alpha$-Catoni-type CS} & \multicolumn{3}{c}{$\alpha$-Catoni-type CS of WR}                                              \\
      \cmidrule(lr){2-4} \cmidrule(lr){5-7}
      $\delta$ & $t=100$                                              & $t=1000$                                          & $t=10000$ & $t=100$ & $t=1000$ & $t=10000$ \\
      \midrule
      0.2      & 0.80742                                              & 0.40365                                           & 0.21584   & 1.12657 & 0.57772  & 0.31491   \\
      0.1      & 0.90313                                              & 0.44475                                           & 0.23620   & 1.26359 & 0.63759  & 0.34483   \\
      0.05     & 0.99164                                              & 0.48093                                           & 0.25378   & 1.38955 & 0.68974  & 0.37059   \\
      0.02     & 1.09822                                              & 0.52301                                           & 0.27411   & 1.54573 & 0.75137  & 0.40032   \\
      0.01     & 1.17616                                              & 0.55114                                           & 0.28791   & 1.66006 & 0.79273  & 0.42069   \\
      0.005    & 1.24275                                              & 0.57863                                           & 0.30049   & 1.75598 & 0.83180  & 0.43922   \\
      0.002    & 1.33019                                              & 0.61138                                           & 0.31623   & 1.88802 & 0.88060  & 0.46188   \\
      0.001    & 1.39427                                              & 0.63528                                           & 0.32692   & 1.98353 & 0.91671  & 0.47799   \\
      \bottomrule
    \end{tabular}
  }
  \caption{Mean interval widths for the non-stitched CSs under the centered Pareto$(1.8)$ model (upper panel) and the centered $t_2$ model (lower panel), with $\alpha=0.5$. Each entry is the empirical mean over $200$ independent runs for the indicated pair $(t,\delta)$. }
  \label{delta}
\end{table}

\FloatBarrier

For both distributions, the empirical mean interval widths increase only mildly as $\delta$ decreases. This is consistent with the $\mathcal{O}\bigl((\log(1/\delta))^{\alpha/(1+\alpha)}\bigr)$ dependence on the confidence level in the explicit width upper bound, within the range where its sufficient condition holds.
At every reported value of $t$ and $\delta$, the proposed procedure yields smaller empirical mean interval widths than the Wang--Ramdas procedure. The gain is especially pronounced when $t$ is moderate and remains visible even at $t=10000$.

\begin{figure}[htbp]
  \centering
  \captionsetup{size=footnotesize, skip=4pt}
  \begin{subfigure}[b]{0.495\textwidth}
    \includegraphics[width=\textwidth]{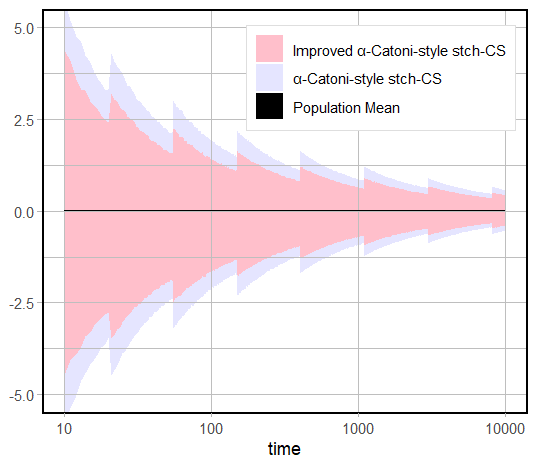}
    \caption*{Centered Pareto$(1.8)$ model}
  \end{subfigure}
  \hfill
  \begin{subfigure}[b]{0.495\textwidth}
    \includegraphics[width=\textwidth]{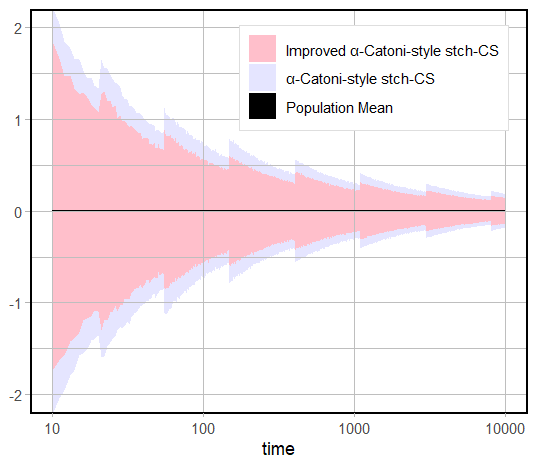}
    \caption*{Centered $t_2$ model}
  \end{subfigure}
  \caption{Stitched procedures at level $\delta=0.05$ with $\alpha=0.5$. The lower and upper boundaries of each band are the empirical means of the corresponding interval endpoints over $200$ independent runs. Their vertical separation is the empirical mean interval width. The visible jumps in the band boundaries reflect transitions between time blocks.  In both models, the proposed stitched CS again has a smaller empirical mean interval width than its Wang--Ramdas counterpart over the displayed range.}
  \label{fig:stitched-width}
\end{figure}

Finally, we compare the stitched procedures. We again fix $\delta=0.05$ and simulate i.i.d. data separately from the same two models. The corresponding blockwise tuning sequences are
\begin{equation*}
  \Theta_{k}=\bigg(\frac{2\log (2 / \delta _{k})}{ \alpha \, C_{\alpha} \, \nu _\alpha (1+u^{-\alpha }) \, \mathrm{e}^{k+1}}\bigg)^{\frac{1}{1+\alpha } }
  \qquad\text{and}\qquad
  \Theta_{k}^{\rm (WR)}=\bigg(\frac{2(1+\alpha)\log (2 / \delta _{k})}{ 5\,\alpha \, \nu _\alpha \, \mathrm{e}^{k+1}}\bigg)^{\frac{1}{1+\alpha } },
\end{equation*}
where the proposed procedure uses the benchmark value $u=u_0$, as in Remark~\ref{choose of omega in SM}. Here $\{\Theta_k\}$ and $\{\Theta_k^{\rm (WR)}\}$ are the blockwise deterministic sequences defined in \eqref{eq:stitched-tuning} and \eqref{eq:wr-stitched-tuning}, respectively. Figure \ref{fig:stitched-width} reports, for each $t$, the empirical mean lower and upper endpoints of the stitched intervals over $200$ independent runs. Their vertical separation is the empirical mean stitched interval width.

\FloatBarrier

\FloatBarrier

The empirical mean endpoint curves are piecewise smooth and exhibit jumps at the boundaries of the time blocks, where the empirical mean interval widths increase. This behavior is intrinsic to the stitching construction and should not be interpreted as a numerical artifact. Between consecutive boundaries, the empirical mean widths continue to decrease with $t$. In both models, the proposed stitched CS again has a smaller empirical mean interval width than its Wang--Ramdas counterpart over the displayed range.

Turning to empirical validity, no cumulative miscoverage event is observed for either stitched procedure in the present experiments: over $4000$ independent runs, neither method exhibits a coverage failure under either model over the displayed monitoring horizon. We therefore omit the corresponding cumulative-miscoverage plots. In this simulation setting, the stitched procedures are empirically more conservative than their non-stitched counterparts. Such empirical conservativeness is consistent with the time-uniform coverage guarantee: the nominal error level is an upper bound on the miscoverage probability, not a value that must be attained.

For reference, Table~\ref{table:sufficient-sample-sizes} gives the sufficient thresholds for time-uniform width control in the present regime. Under the benchmark choice $(u,\lambda)=(u_0,\lambda_0)$, the threshold for the proposed bound is $t=76$ when $\alpha=0.5$ and $\delta=0.05$. This threshold concerns the theoretical time-uniform width bound; it does not identify an empirical onset time.

%% file: section/Appendix.tex
\section{Auxiliary proofs for the width results}
\label{appendix:proof of thm WR stitching}

\subsection{Proofs of almost sure shrinkage}\label{app:almost-sure-shrinkage}~\\

We prove the almost sure shrinkage statements in Remarks~\ref{rem:width-conditions} and~\ref{joint dependence} using the first Borel--Cantelli lemma. The first result concerns a single stitched CS at a fixed error level, while the second allows the error level to vary along a deterministic sequence of times under the additional summability assumption.

\begin{proof}[Proof of almost sure shrinkage in Remark~\ref{rem:width-conditions}]
    We fix $\delta\in(0,1)$ and take all other parameters as in Theorem~\ref{stitching CS}. For $k\geqslant1$ and $j\in\{0,1\}$, let $I_{k,j}=[t_k\mathrm{e}^{j/2},t_k\mathrm{e}^{(j+1)/2})$. If $I_{k,j}\cap\mathcal T(H)=\varnothing$, set $E_{k,j}=\varnothing$; otherwise, take $E_{k,j}$ to be the bound-crossing event in \eqref{eq:stitched-failure-events}. By \eqref{eq:stitched-error-budget}, we have
\[
\sum_{k\geqslant1}\sum_{j=0}^1\mathsf P(E_{k,j})
\leqslant\frac2\kappa \, \sum_{k\geqslant1}\delta_k=\delta<\infty.
\]

    The first Borel--Cantelli lemma implies that almost surely only finitely many of these events occur. Since $\tau(H)<\infty$, \eqref{eq:stitched-admissible-times} shows that almost surely the width bound \eqref{stitching result} holds for all sufficiently large integers $t$. The asserted rate and convergence to zero follow from the relation 
    \[
            B_\alpha(t,\delta)
            = \mathcal{O}\mkern-3mu\left(\left(\frac{\log\log t}{t}\right)^{\frac{\alpha}{1+\alpha}}\right) \ 
            \longrightarrow \ 0.
        \qedhere\]
\end{proof}

\begin{proof}[Proof of almost sure shrinkage in Remark~\ref{joint dependence}]
    With the notation in \eqref{eq:joint-regime-notation}, $D_n=o(t_n)$ implies that the inequality in \eqref{eq:stitched-threshold}, with $h=H$ and $(t,\delta)=(t_n,\delta_n)$, holds for all sufficiently large $n$. By Theorem~\ref{stitching CS}, the probability that the width at $(t_n,\delta_n)$ exceeds its stated bound is at most $\delta_n$. The assumption $\Sigma_{n=1}^{\infty}\,\delta_n<\infty$ and the first Borel--Cantelli lemma therefore imply that this bound holds eventually almost surely. Since $B_\alpha(t_n,\delta_n)\to0$, the widths converge to zero almost surely.
\end{proof}

\subsection{Proof of Theorem~\ref{stitching CS of Ramdas}}\label{app:wr-stitched}~\\

The proof of Theorem~\ref{stitching CS of Ramdas} follows the same stitching argument as that of Theorem~\ref{stitching CS}. It uses \cite[Theorem 12]{CS} for the validity of the underlying $\alpha$-Catoni-type CS, the conditional shifted-moment bound \eqref{eq:shifted-moment-bound} for the time-uniform width analysis, and Lemma~\ref{stitching lemma} with $L_\alpha$ and $\Theta_k$ replaced by $L_{\alpha}^{\rm (WR)}$ and $\Theta_k^{\rm (WR)}$ from \eqref{eq:wr-stitched-tuning}, respectively. Throughout the proof, $B_\alpha(t,\delta)$ is as defined in \eqref{B alpha definition}.

\begin{proof}[Proof of Theorem~\ref{stitching CS of Ramdas}]
    We first establish coverage. For each $k\in\mathbb N_0$, \cite[Theorem 12]{CS} gives a $(1-\delta_k)$-CS with constant weight $\Theta_k^{\rm(WR)}$. Applying Lemma~\ref{Stitching method} with the same partition and error allocation as in the proof of Theorem~\ref{stitching CS} gives the claimed coverage for all $t\in\mathbb N$.

    For the time-uniform width bound, we follow the proof of Theorem~\ref{stitching CS} with the influence function $\psi_\alpha^{\mathrm{(C)}}$, the coefficient $1/(1+\alpha)$ and $\lambda=\alpha$, and use the conditional shifted-moment bound \eqref{eq:shifted-moment-bound} at $u=u_0$. For each fixed comparison point $m$, the supermartingale verification in the proof of Lemma~\ref{catoni supermartingales} then applies to the analogue of \eqref{eq:compensated-supermartingales}, with compensator per observation expressed as follows:
    \[
        \frac{(\Theta_k^{\rm(WR)})^{1+\alpha}}{1+\alpha}
        \bigl(4\nu_\alpha+|\mu-m|^{1+\alpha}\bigr).
    \]
    Together with the CS boundary, its central-moment term gives $L_\alpha^{\rm(WR)}=5\nu_\alpha/(1+\alpha)$.

    Fix $k\in\mathbb N$ and $j\in\{0,1\}$, and use $\varepsilon_k$ and $L(\varepsilon_k,\delta_k)$ from \eqref{eq:stitched-proof-errors}. Write $g_t^{\rm(WR)},\allowbreak B_t^{\pm,\rm(WR)},\allowbreak K_t^{\rm(WR)},\allowbreak J_t^{\rm(WR)}$ for the analogous quantities in \eqref{eq:width-score}, \eqref{eq:width-envelopes} and \eqref{eq:width-root-coefficients} with the above influence function and compensator, $\theta_i=\Theta_k^{\rm(WR)}$, $\beta=\varepsilon_k$ and coverage error $\delta_k$.

    Recall that with $I_{k,j}=[t_k\mathrm{e}^{j/2},t_k\mathrm{e}^{(j+1)/2})$,  if  $I_{k,j}\cap\mathcal T(H_{\mathrm{WR}})=\varnothing$, we set $E_{k,j}^{\rm(WR)}=\varnothing$. Otherwise, choose the deterministic real number $s$ as in that proof, now with $H_{\mathrm{WR}},L_\alpha^{\rm(WR)},\Theta_k^{\rm(WR)}$:
    \[
        s=\max\left\{t_k\,\mathrm{e}^{{j}/{2}}, \ 
        \frac{L(\varepsilon_k,\delta_k)}
        {H_{\mathrm{WR}}-L_\alpha^{\rm(WR)}(\Theta_k^{\rm(WR)})^{1+\alpha}}\right\}.
    \]

    The denominator is positive, and $t_k\mathrm{e}^{j/2}\leqslant s\leqslant t<t_k\mathrm{e}^{(j+1)/2}$ for every $t\in I_{k,j}\cap\mathcal T(H_{\mathrm{WR}})$. At this real $s$, the deterministic coefficient formulas give the expressions 
    \[
        K_s^{\rm(WR)}=\frac{(\Theta_k^{\rm(WR)})^\alpha}{1+\alpha} \quad \mbox{and} \quad
        J_s^{\rm(WR)}=L_\alpha^{\rm(WR)}(\Theta_k^{\rm(WR)})^\alpha+
        \frac{L(\varepsilon_k,\delta_k)}{s \, \Theta_k^{\rm(WR)}}.
    \]

    Hence
    \[
        \Theta_k^{\rm(WR)}J_s^{\rm(WR)}
        \leqslant H_{\mathrm{WR}}=\frac{\alpha^{\frac{1}{\alpha}}}{1+\alpha}.
    \]

    Since $(1+\alpha)^\alpha H_{\mathrm{WR}}^\alpha=\alpha$, the deterministic comparison following \eqref{parameter condition} applies with $\lambda=\alpha$ to the fixed points  $\rho^{\rm(WR)}$ and $\varrho^{\rm(WR)}$, where 
    \[
        \rho^{\rm(WR)}=\mu+(1+\alpha)J_s^{\rm(WR)} \quad \mbox{and} \quad 
        \varrho^{\rm(WR)}=\mu-(1+\alpha)J_s^{\rm(WR)}.
    \]

    Indeed, $K_t^{\rm(WR)}=K_s^{\rm(WR)}$ and $J_t^{\rm(WR)}\leqslant J_s^{\rm(WR)}$ for $t\geqslant s$, so for every such integer $t$, we have 
    \[
        B_t^{+,\rm(WR)}(\rho^{\rm(WR)})
        =-B_t^{-,\rm(WR)}(\varrho^{\rm(WR)})
        \leqslant-\frac{\nu_\alpha}{1+\alpha}\,t\,(\Theta_k^{\rm(WR)})^{1+\alpha}
        -\log\frac2{\delta_k}.
    \]

    Define $E_{k,j}^{\rm(WR)}$ as in \eqref{eq:stitched-failure-events}, with $g_t,B_t^\pm,\rho,\varrho$ replaced by their WR counterparts. Applying Ville's inequality to the corresponding two supermartingales at the fixed points $\rho^{\rm(WR)}$ and $\varrho^{\rm(WR)}$, with boundary $2/\varepsilon_k$ for each, gives $\mathsf P(E_{k,j}^{\rm(WR)})\leqslant\varepsilon_k$. These processes use all observations starting from $X_1$. The same strict monotone-exclusion argument therefore gives, on $(E_{k,j}^{\rm(WR)})^c$,
    \[
        \mathrm{CI}_t^{\mathrm{stch(WR)}}\subset
        (\varrho^{\rm(WR)},\rho^{\rm(WR)})
        \quad\text{for all } \ t\in I_{k,j}\cap\mathcal T(H_{\mathrm{WR}}).
    \]

    Using $L(\varepsilon_k,\delta_k)\leqslant2\log(2/\delta_k)$ and Lemma~\ref{stitching lemma} with $L_\alpha$ replaced by $L_\alpha^{\rm(WR)}$, on this same event we obtain
    \[
        \begin{aligned}
                         \left|\mathrm{CI}_t^{\mathrm{stch(WR)}}\right|
            \leqslant{} &  2(1+\alpha)J_s^{\rm(WR)}                            \\
            \leqslant{} & 2(1+\alpha)
            \frac{L_\alpha^{\rm(WR)}s\,(\Theta_k^{\rm(WR)})^{1+\alpha}+2\log(2/\delta_k)}
            {s\,\Theta_k^{\rm(WR)}}                                        \\
            \leqslant{} & 2A_\alpha^{\rm(WR)}B_\alpha(t,\delta), 
        \end{aligned}
    \]
    simultaneously for all such $t$.

    Finally, the same error budget for the time subintervals as in \eqref{eq:stitched-error-budget} gives
    \[
        \mathsf P\left(\bigcup_{k\geqslant1}\bigcup_{j=0}^1E_{k,j}^{\rm(WR)}\right)
        \leqslant\frac2\kappa\sum_{k\geqslant1}\delta_k=\delta.
    \]

    The width bound thus holds simultaneously on $\mathcal T(H_{\mathrm{WR}})$ with probability at least $1-\delta$. By \eqref{eq:stitched-admissible-times}, this is precisely the claimed bound for every integer $t\geqslant\tau(H_{\mathrm{WR}})$.
\end{proof}